\documentclass[reqno]{amsart}
\usepackage{amsmath,amssymb,stmaryrd}
\begin{document}
\newtheorem{thm}{Theorem}[section]
\newtheorem{lem}{Lemma}[section]
\newtheorem{prop}{Proposition}[section]
\newtheorem{coro}{Corollary}[section]
\newtheorem{defi}{Definition}[section]
\def \eps {\epsilon}
\def \P {{\mathcal P}}
\def \Porth {{\mathcal P}^{\bot}}
\def \Q {{\mathcal Q}}
\def \Qorth {{\mathcal Q}^{\bot}}
\def \rn {\mathbb{R}^n}
\def \rnm {\mathbb{R}^n_-}
\def \rr {\mathbb{R}}
\def \nn {\mathbb{N}}
\def \tu {\tilde{u}}
\def \we {w_\eps}
\def \ze {z_\eps}
\def \he {h_\eps}
\def \xune {x_{1,\eps}}
\def \pe {p_\eps}
\def \tge {\tilde{g}_\eps}
\def \lae {\lambda_\eps}
\def \xe {x_\eps}
\def \ye {y_\eps}
\def \fe {f_\eps}
\def \elle {l_\eps}
\def \re {\rho_\eps}
\def \ae {a_\eps}
\def \ue {u_\eps}
\def \me {\mu_\eps}
\def \ke {k_\eps}
\def \be {\beta_\eps}
\def \ve {v_\eps}
\def \Ve {V_\eps}
\def \bwe {\bar{w}_\eps}
\def \bw {\bar{w}}
\def \bge {\bar{g}_\eps}
\def \Re {R_\eps}
\def \tv {\tilde{v}}
\def \tve {\tilde{v}_\eps}
\def \tue {\tilde{u}_\eps}
\def \txe {\tilde{x}_\eps}
\def \tye {\tilde{y}_\eps}
\def \tze {\tilde{z}_\eps}
\def \bze {\bar{z}_\eps}
\def \crit {2^{\star}}
\def \huno {H_{1,0}^2(\Omega)}
\def \hunrn {H_{1,0}^2(\rn)}
\def \hunrnm {H_{1,0}^2(\rnm)}
\def \ds {\displaystyle}
\def \beq {\begin{eqnarray*}}
\def \eeq {\end{eqnarray*}}
\def \beqn {\begin{eqnarray}}
\def \eeqn {\end{eqnarray}}
\def \bequa {\begin{equation}}
\def \eequa {\end{equation}}

\title[Best Hardy-Sobolev constant]{Elliptic equations with critical
growth and a large set of boundary singularities}

\author{N. Ghoussoub}
\address{Nassif Ghoussoub, Department of Mathematics, University of
British Columbia, Vancouver, Canada}
\thanks{Research partially supported by the Natural Sciences and
Engineering Research Council of Canada.The first named author gratefully
acknowledges the hospitality and support of the Universit\'e de Nice where
this work was initiated.}
\email{nassif@math.ubc.ca}
\author{F. Robert}
\address{Fr\'ed\'eric Robert, Laboratoire J.A.Dieudonn\'e, Universit\'e de
Nice Sophia-Antipolis,
Parc Valrose, 06108 Nice cedex 2, France}
\email{frobert@math.unice.fr}
\thanks{The second named author gratefully acknowledges the hospitality
and support of the University of British Columbia.}
\date{August 17th 2005}
\begin{abstract} We solve variationally certain equations of stellar
dynamics of the form $-\sum_i\partial_{ii} u(x) =\frac{|u|^{p-2}u(x)}{{\rm
dist} (x,{\mathcal A} )^s}$  in a domain $\Omega$ of $\rn$,  where ${\mathcal A} $  is a 
proper linear subspace
of $\rn$. Existence problems are related to the question of 
attainability of the
best constant in the following recent inequality  of 
Badiale-Tarantello \cite{BaTa}:
$$0<\mu_{s,\P}(\Omega)=\inf\left\{\int_{\Omega}|\nabla u|^2\, dx;\;  u\in
\huno \hbox{ and }\int_{\Omega}\frac{|u(x)|^{\crit(s)}}{|\pi(x)|^s}\, 
dx=1\right\}$$
where $0<s<2$, $\crit(s)=\frac{2(n-s)}{n-2}$ and where $\pi$ is the orthogonal projection on a linear space $\P$, where $\hbox{dim}_{\rr}\P \geq 2$. We investigate
this question and how it depends on  the relative position of the subspace
$\Porth$, the orthogonal of $\P$, with respect to the domain $\Omega$ as well as on the curvature of
the boundary $\partial\Omega$ at its points of intersection with $\Porth $.
\end{abstract}
\maketitle

\section{Introduction} Let $\Omega$ be a smooth domain of $\rn$, where
$n\geq 3$, and  denote by $\huno$ the completion of $C_c^\infty(\Omega)$,
the set of smooth functions compactly supported in $\Omega$, for the norm
$\Vert u\Vert_{\huno}=\sqrt{\int_{\Omega}|\nabla u|^2\, dx}.$
%For any subset $\Porth $ of $\rn$, we denote by $d(\cdot,\Porth )$ the Euclidean
%distance to $\Porth $. 
In \cite{BaTa}, Badiale and Tarantello proved that if
$\P $ is a linear subspace of $\rn$ such that $2\leq \hbox{dim}_{\rr}\P \leq n$, then there exists $C>0$ such that for all $u\in \hunrn$,
\bequa\label{ineq:HS:rn}
\left(\int_{\rn}\frac{|u|^{\crit}}{|\pi(x)|^s}\,
dx\right)^{\frac{2}{\crit}}\leq C \int_{\rn}|\nabla u|^2\, dx,
\eequa
where here $\crit=\frac{2(n-s)}{n-2}$, $s\in (0,2)$ and $\pi$ is the orthogonal projection on $\P$ with respect to the Euclidean structure.  Define
\bequa\label{def:mus0}
\mu_{s,\P }(\Omega)=\inf\left\{\frac{\int_{\Omega}|\nabla u|^2\,
dx}{\left(\int_{\Omega}\frac{|u|^{\crit}}{|\pi(x)|^s}\,
dx\right)^{\frac{2}{\crit}}}; u\in \huno\setminus\{0\}\right\}
\eequa
and note that  (\ref{ineq:HS:rn}) and (\ref{def:mus0}) give that for all
smooth domain $\Omega\subset \rn$, we have
\bequa\label{ineq:mus:1}
\mu_{s,\P }(\Omega)\geq \mu_{s,\P }(\rn)>0.
\eequa
   In this article, we adress the question of the value of  the best
constant $\mu_{s,\P }(\Omega)$ as well as the issue of its attainability.
As we will see, both questions are closely related to the relative
positions of $\Porth $ and $\Omega$, and to the geometry of the boundary
$\partial \Omega$ on the points of $\Porth \cap \partial \Omega$.

\medskip\noindent The case when $s=0$ (i.e., the non-singular case) is the
well known Sobolev inequality. In this situation the infimum
$\mu_{s,\P }(\Omega)=\mu_{0,\P }(\rn)$ is not attained unless $\Omega$ is
essentially the whole of $\rn$.

\medskip\noindent The case $s\in (0,2)$ and $\hbox{dim}_{\rr}\P =n$ (that
is $\P =\rn$) was tackled in \cite{gk}, \cite{gr1}, 
\cite{gr2}. It was proved that when $0\in\partial\Omega$, the
infimum in \eqref{def:mus0} is then attained as soon as the mean curvature
of $\partial\Omega$ (oriented with outward pointing normal vectors) at
$0$ is negative. The proof of this result required refined asymptotics
for blown-up solutions of associated second order elliptic equations, 
the difficult case being when these solutions develop a ``bubble" 
located precisely at the point $0$. However, the bubble inherits the
symmetry properties of the problem, and this allowed us to show in 
\cite{gr1} that mean curvature
conditions --as opposed to sectional curvature-- suffice to eliminate the
possibility of a bubbling-off phenomenon.

%\medskip\noindent

In the present paper, we tackle the case of a larger affine subspace of
singularities ($1\leq \hbox{dim}_{\rr}\P \leq n-1$) and in particular when
$\Porth$ contains at least a line.  The situation here closely depends on the
relative positions of $\Porth $ and $\Omega$, the most interesting case being
when the subspace $\Porth $ does not touch the domain $\Omega$ but does touch
its boundary (i.e., when  $\Porth \cap\Omega=\emptyset$ and $\Porth \cap\partial
\Omega\neq\emptyset$). A large part of the analysis is similar to what we
have done in \cite{gr1,gr2} for the case of a single point of singularity
on the boundary of $\Omega$.  However, a new set of difficulties arise  in
this situation: for one,  the centers of the appearing bubbles are not
bound to any particular location and may appear anywhere on
$\partial\Omega$. They do eventually  converge to a point in
$\Porth \cap\partial\Omega\neq \emptyset$, and an important  new issue becomes
the precise control of the distance between the center of the bubble and
this limiting point.

Another new problem related to this setting  is the lack of symmetry of
the bubble. As described by the next proposition, we do show  that it 
enjoys the best symmetry possible in  the $\P$-direction. 
Here and in the sequel,
$\Delta=-\sum_i\partial_{ii}$ will denote the Laplacian with minus sign
convention and $\rnm=\{x\in\rnm/\, x_1<0\}$.

\begin{prop}\label{prop:sym} Let $\pi$ be the projection on a linear
subspace $\Q$ of $\rn$ such that $2\leq \hbox{dim}_{\rr}\Q$ 
and $\Qorth\subset \partial\rnm$.
Assume $s\in (0,2)$ and consider $u\in C^2(\rnm)\cap
C^1(\overline{\rnm})$ such that
\bequa\label{sys:sym}
\left\{\begin{array}{ll}
\Delta u=\frac{u^{\crit (s)-1}}{|\pi (x)|^s}& \hbox{ in }{\rnm}\\
u>0&\hbox{ in }\rnm\\
u=0&\hbox{ on }\partial\rnm.
\end{array}\right.
\eequa
and for some $C>0$,
  \bequa\label{sys:sym.bound}
\hbox{$ u(x)\leq
C(1+|x|)^{1-n}$ for all $x\in\rnm$.}
\eequa
  Then  there
exists $v\in C^2(\rr_-^\star\times \rr\times \Qorth)\cap C^1(\rr_-\times
\rr\times \Qorth)$ such
that for all $z\in \Qorth$, and all $x_1<0$ and $y\in\rn$ with 
$(x_1,y)\in \Q$, we
have that
$u(x_1,y,z)=v(x_1,|y|, z)$.
\end{prop}
But this is not sufficient since the  behavior of the bubble in the
$\P$-direction and the ${\Porth }$-direction cannot often be related.
Overcoming these difficulties, we prove the following theorem. In the 
sequel, $T_x\partial\Omega$ denotes the linear tangent space of 
$\partial\Omega$ at the point $x$.
\begin{thm}\label{th:intro}
Let $\Omega$ be a smooth bounded oriented domain of $\rn$, $n\geq 3$, and
let $\P$ be a linear subspace of $\rn$ such that $2\leq
\hbox{dim}_{\rr}\P $. Assume $s\in (0,2)$.

\begin{trivlist}
\item {\rm (A)}\,  If $\Porth  \cap\Omega\neq\emptyset$, then
$\mu_{s,\P }(\Omega)=\mu_{s,\P }(\rn)$ and the infimum in (\ref{def:mus0})
is not achieved.
\item \item {\rm (B)}\,  If $\Porth  \cap\overline{\Omega}=\emptyset$, then the
infimum  in (\ref{def:mus0}) is achieved.
\item  \item {\rm (C)}\, If $\Porth  \cap\Omega=\emptyset$ and $\Porth 
\cap\partial\Omega\neq\emptyset$, then  the infimum  in (\ref{def:mus0})
is achieved and  the set of minimizers is pre-compact in $\huno$, provided
that at any point $x\in \Porth \cap\partial\Omega$ the principal curvatures of
$\partial\Omega$ at $x$ are non-positive, but do not all vanish. \\
Moreover, at those points $x\in \Porth  \cap\partial\Omega$ where 
$\P\cap T_x\partial\Omega$ and
${\Porth }$ are orthogonal with respect to the second fundamental form of
$\partial \Omega$ at $x$, it suffices that the mean curvature vector of
$\partial\Omega\cap (\Porth +(T_x\partial\Omega)^{\bot}) $ at $x$ be 
null, while the mean curvature of
$\partial\Omega$ at $x$ is negative.
   \end{trivlist}
\end{thm}
The second part in (C) makes connection with the case where $\P=\rn$ (i.e., $\Porth=\{0\} $ studied in \cite{gr1}.
Then the negativity of the mean curvature of $\partial \Omega$ at that
point is sufficient for  $\mu_{s,\P }(\Omega)$ to be attained. On the
opposite end, one may ask  what happens in the case
$\hbox{dim}_{\rr}\P \in \{0,1\}$. In the case when $\P =\{0\}$, inequality
\eqref{ineq:HS:rn} is  clearly irrelevant, however the case
$\hbox{dim}_{\rr}\P =1$ presents some  interest, and this is the object
of the following proposition:

\begin{prop}\label{prop:n-1}
Let $\Omega$ be a smooth bounded oriented domain of $\rn$, $n\geq 3$, and
let $\P$ be a linear subspace of $\rn$ such that $
\hbox{dim}_{\rr}\P =1$. Assume $s\in (0,2)$.

\begin{trivlist}
\item {\rm (A)}\,  If $\Porth  \cap\Omega\neq\emptyset$, then the infimum in
(\ref{def:mus0}) is not achieved.

%Moreover,
%$\mu_{s,\P }(\Omega)=\mu_{s,\P }(\rn)>0$ if $s\in (0,1)$,  and
%$\mu_{s,\P }(\Omega)=0$ whenever $s\in [1,2)$.

   \item  \item {\rm (B)}\,  If $\Porth  \cap\overline{\Omega}=\emptyset$, then
the infimum $\mu_{s,\P }(\Omega)$ in (\ref{def:mus0}) is positive and is
achieved.

\item  \item {\rm (C)}\, If $\Porth  \cap\Omega=\emptyset$ while $\Porth 
\cap\partial\Omega\neq\emptyset$, then $\mu_{s,\P }(\Omega)>0$ and the
infimum is not achieved.
\end{trivlist}
\end{prop}
%\medskip\noindent

Actually, when dealing with case (C) of Theorem \ref{th:intro} and
Proposition \ref{prop:n-1}, the crucial point is to have negative
principal curvatures at each point of $\Porth \cap \partial\Omega$. But the
fact that $\Porth $ only touches $\overline{\Omega}$ at its boundary  means that the
principal curvatures  in the ${\Porth }-$direction are all nonnegative at these
points --at least for those where ${\Porth }$ and $\P\cap 
T_x\partial\Omega$ are orthogonal for the
fundamental form of $\partial \Omega$: therefore,  for
$\mu_{s,\P }(\Omega)$ to be achieved, one needs the negativity of the
principal curvatures in some of the orthogonal directions, which is
obviously impossible if $\Porth $ is $(n-1)-$dimensional and therefore the best
constant is never achieved in this case. This means that the dimension
restriction on the linear subspace in Theorem \ref{th:intro} is optimal.
As a consequence of the techniques developed for the proof of Theorem
\ref{th:intro}, we get the following corollary.
\begin{coro}\label{th:eq}
Let $\Omega$ be a smooth bounded oriented domain of $\rn$ and let $\pi$ be
the orthogonal projection onto a linear vector subspace $\Q\subset\rn$
such that $2\leq \hbox{dim}_{\rr}\Q$. We assume that $\Qorth
\cap\Omega=\emptyset$ and $\Qorth \cap\partial\Omega\neq\emptyset$. Assume
$s\in (0, 2)$ and consider $a\in C^1(\overline{\Omega})$ such that the
operator $\Delta+a$ is coercive on $\Omega$. Then there exists a solution
$u\in \huno\cap C^1(\overline{\Omega})$ for
$$\left\{\begin{array}{ll}
\Delta u+au=\frac{u^{\crit-1}}{|\pi (x)|^s}& \hbox{ in }{\mathcal
D}'(\Omega)\\

u>0&\hbox{ in }\Omega\\

u=0&\hbox{ on }\partial\Omega

\end{array}\right.$$
   provided that at any point $x\in \Qorth \cap\partial\Omega$   the
principal curvatures of $\partial\Omega$ at $x$ are non-positive, but do
not all vanish.\\
    Moreover, at those points $x\in \partial \Omega \cap \Qorth$ where
$\Qorth$ and $\Q\cap T_x\partial\Omega $ are orthogonal with respect 
to the second fundamental
form of $\partial \Omega$  at $x$, it suffices to assume that the mean
curvature vector of 
$\partial\Omega\cap(\Qorth+(T_x\partial\Omega)^{\bot})$ at $x$ 
is null, while the mean
curvature of $\partial\Omega$ at $x$ is negative.
\end{coro}

Related references for best constant problems in Sobolev inequalities are Druet \cite{d1}, Hebey-Vaugon \cite{hv1,hv2} and Egnell \cite{eg}. Concerning asymptotics for blown-up sequences of solutions to elliptic equations, we also refer to Atkinson-Peletier \cite{ap}, Br\'ezis-Peletier \cite{bp}, Han \cite{han}, Druet \cite{d2}, Druet-Hebey \cite{dh2}, Druet-Hebey-Robert \cite{dhr} and Schoen-Zhang \cite{sz}.

The rest of the paper is devoted to the proof of these results. As
mentioned above, a significant part of the analysis was developed in
\cite{gr1,gr2} for the case of a unique singular point at the boundary,
and to which we shall refer frequently.  On the other hand, we shall give
all the details relating to the new difficulties arising in this new
setting of large set of singularities. The paper is organized as follows.
In section 1, we deal with points (A) and (B) of theorem 
\ref{th:intro} and prove a symmetry result. Sections 2
to 5 are devoted to the proof of point (C) of Theorem \ref{th:intro} 
which is much more intricate,
as it will  require the full range of modern techniques for blow-up
analysis and strong pointwise estimates for  minimizers of the
subcritical functional associated to (\ref{def:mus0}). In section 6, 
we prove Proposition \ref{prop:n-1}, while the appendix in section 7 
provides a
required regularity result for the family of elliptic pde's with
singularities that we are dealing with in this paper. As a last remark, note that all the statements can be straightforwardly adapted to the case when $\P$ is an affine subspace of $\rn$, and not only a linear space.

\section{Partial symmetry of bubbles and Part (A), (B) of Theorem
\ref{th:intro}}
%test-functions estimates}\label{sec:test}
We let $\P$ be a linear subspace of $\rn$ with $2\leq \hbox{dim}_{\rr}\P\leq n-1$. We shall denote by $\pi$ the orthogonal 
projection on
$\P $, and 
\bequa\label{def:mus}
\mu_{s,\P }(\Omega)=\inf\left\{\frac{\int_{\Omega}|\nabla u|^2\,
dx}{\left(\int_{\Omega}\frac{|u|^{\crit}}{|\pi(x)|^s}\,
dx\right)^{\frac{2}{\crit}}}; u\in \huno\setminus\{0\}\right\}
\eequa

\medskip\noindent{\bf Proof of Proposition \ref{prop:sym}.} We first 
prove the partial symmetry property for the positive solutions to the 
limit
equation on
$\rnm$.
   For that, we consider $u\in C^2(\rnm)\cap C^1(\overline{\rnm})$ that
verifies the system (\ref{sys:sym}) while verifying for some $C>0$ the
bound $u(x)\leq \frac{C}{(1+|x|)^{n-1}}$.
We follow the proof of \cite{gr1} to which we refer for details. For 
simplicity, up to a change of coordinates, we write any point 
$x\in\rn$ as $x=(x_1,y,z)$, where $(x_1,y)\in \Q=\rr^{k}$ and 
$z\in\Qorth=\rr^{n-k}$. Therefore $\pi(x)=(x_1,y,0)$. We let 
$\vec{e}_1$ be the first vector of the canonical
basis of $\rn$ and consider  the open ball
$$D:=B_{1/2}\left(-\frac{1}{2}\vec{e}_1\right).$$
We define
\bequa
v(x):=|x|^{2-n}u\left(\vec{e}_1+\frac{x}{|x|^2}\right)
\eequa
for all $x\in \overline{D}\setminus\{0\}$. We extend $v$ by $0$ at
$0$. This is then well-defined and $v\in C^2(D)\cap 
C^1(\overline{D}\setminus\{0\})$. Moreover, $v(x)>0$ for all $x\in D$ 
and
$v(x)=0$ for all $x\in\partial D\setminus\{0\}$. The function $v$ 
verifies the equation
\bequa\label{eq:v}
\Delta
v=\frac{v^{\crit-1}}{|\pi(x+|x|^2\vec{e}_1)|^s}
\eequa
in $D$. Since $v>0$ in $D$, it follows from Hopf's Lemma that
$\frac{\partial v}{\partial \nu}<0$ on $\partial D\setminus\{0\}$.

We prove the symmetry of $u$ by proving a symmetry property of $v$, which
is defined on a ball. Our proof uses the moving plane method. We take
largely inspiration in \cite{gnn} and \cite{cgs}. We let 
$i\in\{2,...,k\}$. For any
$\mu\geq 0$ and $x\in \rn$, we let
$$x_\mu=(x_1,...,2\mu-x_i,...,x_n)\hbox{ and 
}D_\mu=\{x\in D/\, x_\mu\in D\}.$$
It follows from Hopf's Lemma that there exists
$\epsilon_0>0$ such that for any $\mu\in
(\frac{1}{2}-\epsilon_0,\frac{1}{2})$, we have that $D_\mu\neq\emptyset$
and $v(x)\geq v(x_\mu)$ for all $x\in D_\mu$ such that $x_i\leq\mu$. We
let $\mu\geq 0$. We say that $(P_\mu)$ holds if $D_{\mu}\neq \emptyset$
and $v(x)\geq v(x_{\mu})$ for all $x\in D_{\mu}$ such that $x_i\leq\mu$. We let
\bequa\label{def:lambda}
\lambda:=\min\left\{\mu\geq 0;\, (P_{\nu})\hbox{ holds for all }\nu\in
\left(\mu,\frac{1}{2}\right) \right\}.
\eequa

\medskip We claim that $\lambda=0$. Indeed we proceed by 
contradiction and assume that $\lambda>0$. We
then get that $D_{\lambda}\neq\emptyset$ and that $(P_{\lambda})$ holds.
We let $w(x):=v(x)-v(x_{\lambda})$
for all $x\in D_{\lambda}\cap\{x_n<\lambda\}$. Since $(P_{\lambda})$
holds, we have that $w(x)\geq 0$ for all $x\in
D_{\lambda}\cap\{x_i<\lambda\}$. With the equation (\ref{eq:v}) of $v$ and
$(P_{\lambda})$, we get that
\beq
\Delta w&=& \frac{v(x)^{\crit-1}}{|\pi(x+|x|^2\vec{e}_1)|^s}-
\frac{v(x_{\lambda})^{\crit-1}}{|\pi(x_{\lambda}+|x_{\lambda}|^2\vec{e}_1)|^s}\\
&\geq &v(x_{\lambda})^{\crit-1}\left(\frac{1}{|\pi(x+|x|^2\vec{e}_1)|^s}-
\frac{1}{|\pi(x_{\lambda}+|x_{\lambda}|^2\vec{e}_1)|^s}\right)
\eeq
for all $x\in D_{\lambda}\cap\{x_i<\lambda\}$. Since $2\leq i\leq k$, 
we get that the RHS is positive (see \cite{gr1}), and then $\Delta 
w(x)>0$ for all $x\in
D_{\lambda}\cap\{x_i<\lambda \}$. It then follows from Hopf's Lemma and the strong comparison principle that
$$w>0\hbox{ in }D_\lambda\cap\{x_i<\lambda\}\hbox{ and }\frac{\partial
w}{\partial \nu}<0\hbox{ on }D_\lambda\cap\{x_i=\lambda\}.$$
The contradiction then follows from standard arguments, we refer to 
\cite{gr1,gr2} for details. This yields $\lambda=0$.

\medskip\noindent Here goes the final argument.
Since $\lambda=0$, it follows from the definition (\ref{def:lambda}) of
$\lambda$ that $v(x)\geq v(x_1,...,-x_i,...,x_n)$ for all $x\in D$ such that
$x_i\leq 0$. With the same technique, we get the reverse inequality, and
then, we get that $v(x)=v(x_1,...,-x_i,...,x_n)$ for all 
$x=(x',x_n)\in D$. In other words, $v$ is symmetric with respect
to the hyperplane $\{x_i=0\}$. The same analysis holds for any
hyperplane containing $\hbox{Span}\{\vec{e_1},\vec{e_{k+1}},...., 
\vec{e_n}\}$. Coming back to the initial function
$u$, this proves Proposition \ref{prop:sym} and the symmetry property.

\medskip\noindent The object of the following proposition is to deal 
with case (A) of
Theorem \ref{th:intro} that is when $\Porth  \cap\Omega\neq\emptyset$.

\begin{prop}\label{prop:Case1}
Let $\Omega$ be a smooth bounded domain of $\rn$, $n\geq 3$.  Let
$\P \subset\rn$ be a linear vector subspace of $\rn$, where $2\leq
\hbox{dim}_{\rr}\P \leq n-1$. Let $s\in (0,2)$ and assume that $\Porth 
\cap\Omega\neq\emptyset$, then
$\mu_{s,\P }(\Omega)=\mu_{s,\P }(\rn)$
and the infimum $\mu_{s,\P }(\Omega)$ is not achieved.
\end{prop}

\smallskip\noindent{\it Proof:} Fix $x_0\in \Porth  \cap\Omega$, and let
$\delta>0$ such that $B_\delta(x_0)\subset\Omega$. Let $\alpha>0$ and
$u\in  C_c^\infty(\rn)\setminus\{0\}$ such that
$$\frac{\int_{\rn}|\nabla u|^2\,
dx}{\left(\int_{\rn}\frac{|u|^{\crit}}{|\pi(x)|^s}\,
dx\right)^{\frac{2}{\crit}}}\leq \mu_{s,\P }(\rn)+\alpha.$$
For $\eps>0$, we  let
$\ue(x)=\eps^{-\frac{n-2}{2}}u\left(\frac{x-x_0}{\eps}\right)$
for all $x\in\Omega$. As easily checked, $\ue\in C^\infty_c(\Omega)$ 
for $\eps>0$ small and
$$\frac{\int_{\Omega}|\nabla \ue|^2\,
dx}{\left(\int_{\Omega}\frac{|\ue|^{\crit}}{|\pi(x)|^s}\,
dx\right)^{\frac{2}{\crit}}}=\frac{\int_{\rn}|\nabla u|^2\,
dx}{\left(\int_{\rn}\frac{|u|^{\crit}}{\left|\pi(x)\right|^s}\,
dx\right)^{\frac{2}{\crit}}}\leq \mu_{s,\P }(\rn)+\alpha.$$
Here, we have used that $x_0\in\Porth$, that is $\pi(x_0)=0$. Coming back to
the definition (\ref{def:mus}) of $\mu_{s,\P }(\Omega)$ letting $\alpha\to
0$ and using (\ref{ineq:mus:1}), we get that
$\mu_{s,\P }(\Omega)=\mu_{s,\P }(\rn)$.

\smallskip\noindent We claim that $\mu_{s,\P }(\Omega)$ is not achieved.
Indeed,
%we argue by contradiction and
assuming it is achieved by a function $u\in\huno\setminus\{0\}$, we can
assume without loss that $u\geq 0$. Since
$\mu_{s,\P}(\Omega)=\mu_{s,\P}(\rn)$, we get that $\mu_{s,\P}(\rn)$ is
also attained by $u$ which then verifies $\Delta
u=\frac{u^{\crit-1}}{|\pi(x)|^s}$ in ${\mathcal D}'(\rn)$. Since $u\geq
0$, it follows from the regularity results of section \ref{sec:app} 
and the maximum principle that $u>0$ on $\rn\setminus\P$, a
contradiction since $u\in \huno$.\hfill$\Box$\\

The case where $\Porth  \cap\overline{\Omega}=\emptyset$ is dealt with in the
following proposition.
\begin{prop}\label{prop:Case2}
Let $\Omega$ be a smooth bounded domain of $\rn$, $n\geq 3$, and let $\P$
be a linear vector subspace of $\rn$ such that $2\leq
\hbox{dim}_{\rr}\P \leq n-1$. Assume $s\in (0,2)$ and  that $\Porth 
\cap\overline{\Omega}=\emptyset$, then the infimum $\mu_{s,\P }(\Omega)$ is
attained.
\end{prop}

\smallskip\noindent{\it Proof:} Since $\Porth 
\cap\overline{\Omega}=\emptyset$, there exists $c,C>0$ such that $c\leq
|\pi(x)|\leq C$ for all $x\in\Omega$. In particular, since
$\crit<\frac{2n}{n-2}$, we have compactness of the embedding of $\huno$ in
$L^{\crit}(\Omega, |\pi(x)|^{-s})$ and therefore the existence of
minimizers. This ends the proof of the Proposition.\hfill$\Box$

\section{Blow-up analysis, Part I}\label{sec:part1}
Throughout this section, we let $\Omega$ be a smooth bounded domain of
$\rn$, $n\geq 3$, and  $\P$ be a linear vector subspace of $\rn$ such that
$2\leq \hbox{dim}_{\rr}\P\leq n-1$. Let $s\in (0,2)$ and  assume that
\bequa\label{hyp:P:1}
\Porth  \cap\Omega=\emptyset\hbox{ and }\Porth  \cap \partial\Omega\neq\emptyset.
\eequa
Here and in the sequel, we let $\pi$ be the orthogonal projection on 
$\P $.
%, so that $d(x,\Porth )=|\pi(x)|$ for all $x\in\rn$.
%\medskip\noindent We assume in the sequel that
%$\Porth  \cap\Omega=\emptyset\hbox{ and }\Porth  \cap \partial\Omega\neq\emptyset$ which
This is the most intricate case and to which the rest of the paper is
essentially devoted.
\begin{prop}\label{prop:test}
Let $\Omega$ be a smooth bounded domain of $\rn$, $n\geq 3$, and let $\P$ be a linear vector subspace of $\rn$, such that $2\leq
\hbox{dim}_{\rr}\P \leq n-1$. Let  $s\in (0,2)$ and assume that $\Porth 
\cap\Omega=\emptyset$ and $\Porth  \cap \partial\Omega\neq\emptyset$, then
$\mu_{s,\P }(\Omega)\leq \mu_{s,\P }(\rnm).$
\end{prop}
\smallskip\noindent{\it Proof:} Let $x_0\in \Porth  \cap \partial\Omega$. Since
$\Porth  \cap\Omega=\emptyset$, we have that
\bequa\label{incl:porth}
\Porth  \subset T_{x_0}\partial\Omega,
\eequa
where $T_{x_0}\partial\Omega$ is the linear tangent space at $x_0$ of the
smooth manifold $\partial\Omega$. It follows from (\ref{incl:porth}) that
$(T_{x_0}\partial\Omega)^\bot\subset \P $.  We choose a direct
orthonormal basis $(\vec{e}_1,...,\vec{e}_n)$ of $\rn$ such that
\bequa\label{choice:basis}\begin{array}{l}
\vec{e}_1=\vec{n}_{x_0}\hbox{ is the normal outward vector at }x_0\hbox{
of }\partial\Omega\\
(\vec{e}_1,...,\vec{e}_k)\hbox{ is an orthonormal basis of }\P \\
(\vec{e}_{k+1},...,\vec{e}_n)\hbox{ is an orthonormal basis of }\Porth  .
\end{array}\eequa
Here and in what follows, $k=\hbox{dim}_{\rr}\P ,$
so that $2\leq k\leq n-1$. In particular,
$(\vec{e}_{2},...,\vec{e}_n)$ is an orthonormal basis of $T_{x_0}\partial\Omega$. For the rest of this section, we shall be refering to this particular
basis. In particular, we adopt the following notation: we write any
element $x\in\rn $ as $x=(x_1,y,z)$, with $x_1\in\rr$,
$y\in\hbox{span}(\vec{e}_2,...,\vec{e}_k)$ and
$z\in\hbox{span}(\vec{e}_{k+1},...,\vec{e}_n)=\Porth  $.

\medskip\noindent Since $\partial\Omega$ is smooth, there exist $U,V$ open
subsets of $\rn$, such that $0\in U$ and $x_0\in V$, there exists
$\varphi\in C^\infty(U,V)$ and $\varphi_0\in C^\infty(U')$ (where
$U'=\{(y,z)/ \,\exists x_1\in\rr\hbox{ s.t. } (x_1,y,z)\in U\}$) such that
\bequa\begin{array}{ll}
(i) & \varphi: U\to V\hbox{ is a }C^\infty-\hbox{diffeomorphism}\\
(ii) & \varphi(0)=x_0\\
(iii) & \varphi(U\cap\{x_1<0\})=\varphi(U)\cap\Omega\hbox{ and
}\varphi(U\cap\{x_1=0\})=\varphi(U)\cap\partial\Omega.\\
(iv) & \varphi_0(0)=0\hbox{ and }\nabla\varphi_0(0)=0\\
(v) & \varphi(x_1,y,z)=(x_1+\varphi_0(y,z),y,z)+x_0\hbox{ for all
}(x_1,y,z)\in U
\end{array}\label{def:vphi:3}
\eequa
where $D_x\varphi_0$ denotes the differential of $\varphi_0$ at $x$. Let
$\alpha>0$ and $u\in C_c^\infty(\rnm)\setminus\{0\}$ such that
$$\frac{\int_{\rnm}|\nabla u|^2\,
dx}{\left(\int_{\rnm}\frac{|u|^{\crit}}{|\pi(x)|^s}\,
dx\right)^{\frac{2}{\crit}}}\leq \mu_{s,\P }(\rnm)+\alpha.$$
Define
$\ue(x)=\eps^{-\frac{n-2}{2}}u\left(\frac{\varphi^{-1}(x)}{\eps}\right)$
for all $x\in\Omega$ and all $\eps>0$. As easily checked, for $\eps>0$
small enough, we have that $\ue\in C_c^\infty(\Omega)$. Standard
computations yield that
$$\mu_{s,\P }(\Omega)\leq\frac{\int_\Omega |\nabla\ue|^2\,
dx}{\left(\int_\Omega\frac{|\ue|^{\crit}}{|\pi(x)|^s}\,
dx\right)^{\frac{2}{\crit}}}=\frac{\int_{\rnm} |\nabla u|^2\,
dx}{\left(\int_{\rnm}\frac{|u|^{\crit}}{|\pi(x)|^s}\,
dx\right)^{\frac{2}{\crit}}}+o(1)\leq \mu_{s,\P }(\rnm)+\alpha+o(1)$$
where $\lim_{\eps\to 0}o(1)=0$. Letting $\eps\to 0$ and $\alpha\to 0$, we
get the claimed result.
  \hfill$\Box$

\medskip\noindent In order to construct minimizers for
$\mu_{s,\P }(\Omega)$, we consider a subcritical minimization problem for
which we have compactness.  The proof of this result is standard and we
refer to \cite{gr1} for details.

\begin{prop}\label{prop:subcrit} Let $\Omega$ be a smooth bounded domain
of $\rn$, $n\geq 3$, and let $\P$ be a linear vector subspace of $\rn$
such that $2\leq \hbox{dim}_{\rr}\P \leq n-1$. Let $s\in (0,2)$ and assume
that (\ref{hyp:P:1}) holds, then for any $\eps\in (0,\crit-2)$, the
infimum
$$\mu_{s,\P }^\eps(\Omega):=\inf_{u\in\huno\setminus\{0\}}\frac{\int_\Omega
|\nabla u|^2\, dx}{\left(\int_\Omega\frac{|u|^{\crit-\eps}}{|\pi(x)|^s}\,
dx\right)^{\frac{2}{\crit-\eps}}},$$
is achieved by a function $\ue\in\huno$. Moreover, $\ue\in
C^\infty(\overline{\Omega}\setminus\Porth  )$ and can be assumed to satisfy the
system
$$\left\{\begin{array}{ll}
\Delta\ue=\frac{\ue^{\crit-1-\eps}}{|\pi(x)|^s} & \hbox{ in }{\mathcal
D}'(\Omega)\\
\ue>0 & \hbox{ in }\Omega\\
\int_\Omega\frac{\ue^{\crit-\eps}}{|\pi(x)|^s}\,
dx=(\mu_{s,\P }^\eps(\Omega))^{\frac{\crit-\eps}{\crit-2-\eps}}&
\end{array}\right.$$
Moreover, we have that
$\lim_{\eps\to 0}\mu_{s,\P }^\eps(\Omega)=\mu_{s,\P }(\Omega),$
and there exists $u_0\in \huno$ such that, up to a subsequence,
$\ue\rightharpoonup u_0$ weakly in $\huno$ when $\eps\to 0$. If
$u_0\not\equiv 0$, then $\lim_{\eps\to 0}\ue=u_0$ strongly in $\huno$ and
$u_0$ is a minimizer for $\mu_{s,\P }(\Omega)$. In particular,
$\mu_{s,\P }(\Omega)$ is attained.
\end{prop}

  We now start the blow-up analysis for minimizing sequences. 
Actually, we consider a more general case. Here and in
the sequel, we let $\pe\in [0,\crit-2)$ such that
$$\lim_{\eps\to 0}\pe=0.$$
We assume that \eqref{hyp:P:1} holds. We consider a family 
$(\ae)_{\eps>0}\in C^1(\overline{\Omega})$ such that
there exists $\lambda,C>0$ such that
\bequa\label{hyp:ae}
\Vert\ae\Vert_{C^1(\overline{\Omega})}\leq C\hbox{ and
}\int_\Omega(|\nabla\varphi|^2+\ae\varphi^2)\, dx\geq
\lambda\int_\Omega\varphi^2\, dx
\eequa
for all $\eps\to 0$ and all $\varphi\in C^\infty_c(\Omega)$. For any
$\eps>0$, we consider $\ue\in \huno\cap C^2(\overline{\Omega}\setminus \Porth 
)$ a solution to the system
\bequa\label{syst:ue}
\left\{\begin{array}{ll}
\Delta\ue+\ae\ue=\frac{\ue^{\crit-1-\pe}}{|\pi(x)|^s}& \hbox{ in
}{\mathcal D}'(\Omega)\\
\ue>0&\hbox{ in }\Omega
\end{array}\right.
\eequa
We assume that $\ue$ is of minimal energy type, that is
\bequa\label{hyp:nrj:min}
\int_\Omega\frac{|\ue|^{\crit-\pe}}{|\pi(x)|^s}\,
dx=\mu_{s,\P }(\Omega)^{\frac{\crit}{\crit-2}}+o(1)
\eequa
where $\lim_{\eps\to 0}o(1)=0$. We also assume that blow-up occurs, that
is
\bequa\label{hyp:blowup}
\ue\rightharpoonup 0
\eequa
weakly in $\huno$ when $\eps\to 0$. Such a family arises naturally when
$u_0\equiv 0$ in Proposition \ref{prop:subcrit}. It follows from
Proposition \ref{prop:app} of the Appendix that $\ue\in
C^0(\overline{\Omega})$. We let $\xe\in\Omega$ and $\me,\ke>0$ such that
\bequa\label{def:me:xe}
\max_\Omega\ue=\ue(\xe)=\me^{-\frac{n-2}{2}}\hbox{ and
}\ke=\me^{1-\frac{\pe}{\crit-2}}.
\eequa
Our goal in this section is to prove the following:

\begin{prop}\label{prop:sec3:3} Under the above assumption,
  there exists $x_0\in\Porth  \cap\partial\Omega$, a chart $\varphi$ as in
(\ref{def:vphi:3}), there exists $(\bze)_{\eps>0}\in \partial\rnm$ such
that $\lim_{\eps\to 0}\bze=0$ and such that the function
$$\ve(x):=\me^{\frac{n-2}{2}}\ue\circ\varphi(\bze+\ke x)$$
defined for $x\in \frac{U-\bze}{\ke}$ and $\eps>0$ verifies that there
exists $v\in\hunrnm\setminus\{0\}$ such that for any $\eta\in
C_c^\infty(\rn)$,
$\eta\ve\rightharpoonup\eta v\hbox{ in }\hunrnm$ weakly in ${\mathcal 
D}'(\rnm)$ when $\eps\to 0$. The function $v$ verifies that
$$\Delta v=\frac{v^{\crit-1}}{|\pi(x)|^s}\hbox{ in }{\mathcal D}'(\rnm)$$
and
$\int_{\rnm}|\nabla v|^2\,
dx=\mu_{s,\P }(\Omega)^{\frac{\crit}{\crit-2}}=\mu_{s,\P }(\rnm)^{\frac{\crit}{\crit-2}}.$
In addition, $v\in C^{1}(\overline{\rnm})$ and
\bequa\label{cv:v:c1}
\lim_{\eps\to 0}\ve=v\hbox{ in }C^{1}_{loc}(\overline{\rnm}).
\eequa
Moreover,
$$\lim_{\eps\to 0}\me^{\pe}=1.$$
\end{prop}

\smallskip\noindent{\it Proof:} The proof goes in five steps.

\medskip\noindent{\bf Step \ref{sec:part1}.1:} We claim that
\bequa\label{lim:prelim}
\me=o(1)\hbox{ and }\pi(\xe)=O(\ke)
\eequa
when $\eps\to 0$.

Indeed assume  that $\lim_{\eps\to 0}\me\neq 0$, then up to a subsequence,
there exists $C>0$ such that $|\ue(x)|\leq C$ for all $x\in\Omega$ and all
$\eps>0$. Mimicking the proof of the Appendix, we get that there exists
$C>0$ such that $\Vert\ue\Vert_{C^1(\overline{\Omega})}\leq C$. Since
(\ref{hyp:blowup}) holds, it follows from Ascoli's theorem that, up to a
subsequence, $\lim_{\eps\to 0}\ue=0$ in $C^0(\overline{\Omega})$. A
contradiction with (\ref{hyp:nrj:min}). This proves that
$\lim_{\eps\to 0}\me=0.$

To prove the second part of the  claim assume that
\bequa\label{lim:step:1}
\lim_{\eps\to 0}\frac{|\pi(\xe)|}{\ke}=+\infty.
\eequa
For any $\eps>0$, set
\bequa\label{def:be}
\be=|\pi(\xe)|^{\frac{s}{2}}\ue(\xe)^{\frac{2+\pe-\crit}{2}}=|\pi(\xe)|^{\frac{s}{2}}\ke^{\frac{2-s}{2}}.
\eequa
It follows from the definition (\ref{def:be}) of $\be$ and
(\ref{lim:step:1}) that
\bequa\label{ppty:be}
\lim_{\eps\to 0}\be=0,\; \lim_{\eps\to
0}\left(\frac{\be}{\ke}\right)^2=+\infty\hbox{ and }\lim_{\eps\to
0}\left(\frac{\be}{|\pi(\xe)|}\right)^2=0
\eequa
when $\eps\to 0$.

\medskip\noindent{\it Case \ref{sec:part1}.1.1:} Assume first there exists
$\rho>0$ such that
$\frac{d(\xe,\partial\Omega)}{\be}\geq 2\rho$
for all $\eps>0$. For $x\in B_{2\rho}(0)$ and $\eps>0$,  define
$$\ve(x):=\frac{\ue(\xe+\be x)}{\ue(\xe)}.$$
This is well defined since $\xe+\be x\in\Omega$ for all $x\in
B_{2\rho}(0)$. As easily checked, with \eqref{syst:ue}, we have that
$$\Delta\ve+\ke^2\ae(\xe+\be
x)\ve=\frac{\ve^{\crit-1-\pe}}{\left|\frac{\pi(\xe)}{|\pi(\xe)|}+\frac{\be}{|\pi(\xe)|}\pi(x)\right|^s}$$
weakly in $B_{2\rho}(0)$. Since $0\leq \ve(x)\leq \ve(0)=1$ for all
$x\in B_{\rho}(0)$, it follows from standard elliptic theory  and 
\eqref{ppty:be} that there
exists $v\in C^1(B_{2\rho}(0))$ such that
$\ve\to v$
in $C^1_{loc}(B_{2\rho}(0))$ as $\eps\to 0$. In particular,
\bequa\label{lim:v:case1:nonzero}
v(0)=\lim_{\eps\to 0}\ve(0)=1
\eequa
With a change of variables and the definition (\ref{def:be}) of $\be$, we
get that
\beq
&&\int_{\Omega\cap B_{\rho\be}(\xe)}\frac{\ue^{\crit-\pe}}{|\pi(x)|^s}\,
dx=\frac{\ue(\xe)^{\crit-\pe}\be^n}{|\pi(\xe)|^s}\int_{B_{\rho}(0)}\frac{\ve^{\crit-\pe}}{\left|\frac{\pi(\xe)}{|\pi(\xe)|}+\frac{\be}{|\pi(\xe)|}\pi(x)\right|^s}\,dx\\
&&\geq
\left(\frac{\be}{\ke}\right)^{n-2}\me^{-\pe\frac{n-2}{\crit-2}}\int_{B_{\rho}(0)}\frac{\ve^{\crit-\pe}}{\left|\frac{\pi(\xe)}{|\pi(\xe)|}+\frac{\be}{|\pi(\xe)|}\pi(x)\right|^s}\,dx.
\eeq
Using (\ref{hyp:nrj:min}), (\ref{ppty:be}) and passing to the limit
$\eps\to 0$ (note that $\me^{-1}\geq 1$ for $\eps>0$ small), we get that
$\int_{B_{\rho}(0)}v^{\crit}\, dx=0,$
and then $v\equiv 0$. This contradicts  (\ref{lim:v:case1:nonzero}) and
therefore (\ref{lim:step:1}) does not hold, which proves the claim in Case
3.1.1.

\medskip\noindent{\it Case \ref{sec:part1}.1.2:} Now assume that, up to a
subsequence,
$\lim_{\eps\to 0}\frac{d(\xe,\partial\Omega)}{\be}=0.$
We then get a contradiction by a rescaling of $\ue$ as in \cite{gr1}. The
proof uses the techniques of Case \ref{sec:part1}.1.1 and is rather 
similar to \cite{gr1} to which we refer for the details.

\medskip\noindent In both cases, we have obtained a contradiction and Step
3.1 is established.\hfill$\Box$

\medskip\noindent{\bf Step \ref{sec:part1}.2:} Up to a subsequence, we
claim that $x_0$ defined as
\bequa\label{def:x0}
x_0=\lim_{\eps\to 0}\xe
\eequa
belongs to $\Porth  \cap\partial\Omega.$\\

Indeed, it follows from (\ref{lim:prelim}) and \eqref{def:me:xe} that
$\pi(x_0)=0$, that is $x_0\in\Porth  $. Since $x_0\in\overline{\Omega}$, it
follows from (\ref{hyp:P:1}) that $x_0\in \Porth  \cap\partial\Omega$.\\

\medskip\noindent Since (\ref{hyp:P:1}) holds, we have that 
\eqref{incl:porth} holds. We choose a basis as in (\ref{choice:basis})
and we choose a chart $\varphi$ as in (\ref{def:vphi:3}). In particular,
here again, we let $k=\hbox{dim}_{\rr}\P \in \{2,...,n-1\}$.

\medskip\noindent{\bf Step \ref{sec:part1}.3:} Setting
\bequa\label{def:xune:ye:ze}
\xe=\varphi(\xune, \ye,\ze),
\eequa
where $\xune<0$, $\ye\in\hbox{span}(\vec{e}_2,...,\vec{e}_k)$ and
$\ze\in\hbox{span}(\vec{e}_{k+1},...,\vec{e}_n)=\Porth  $, we claim that
\bequa\label{claim:1}
d(\xe,\partial\Omega)=(1+o(1))|\xune|=O(\ke),\; \ye=O(\ke)\hbox{ and
}\varphi_0(0,\ze)=O(\ke),
\eequa
when $\eps\to 0$. Here $\varphi_0$ is as in (\ref{def:vphi:3}).

\smallskip\noindent{\it Proof of the claim:} our first remark is that
\bequa\label{lim:prop32:1}
d(\xe,\partial\Omega)=O(\ke)
\eequa
when $\eps\to 0$. Indeed, since $\Porth  \cap\Omega=\emptyset$, we have that
$\xe-\pi(\xe)\in \Porth  \in \rn\setminus\Omega$. Since $\xe\in\Omega$, there
exists $t_\eps\in (0,1)$ such that
$t_\eps\xe+(1-t_\eps)\cdot(\xe-\pi(\xe))\in \partial\Omega$. Consequently,
$$d(\xe,\partial\Omega)\leq
|\xe-(t_\eps\xe+(1-t_\eps)\cdot(\xe-\pi(\xe)))|=(1-t_\eps)|\pi(\xe)|\leq
|\pi(\xe)|=O(\ke)$$
when $\eps\to 0$. This proves \eqref{lim:prop32:1}.\\

   As in \cite{gr1}, we get that
\bequa\label{lim:prop32:2}
d(\xe,\partial\Omega)=(1+o(1))|\xune|
\eequa
when $\eps\to 0$. We write that
$$\pi(\xe)=\pi(\xune+\varphi_0(\ye,\ze),
\ye,\ze)=(\xune+\varphi_0(\ye,\ze), \ye,0).$$
With (\ref{lim:prelim}) and (\ref{lim:prop32:1}),
we then get that
\bequa\label{eq:lundi}
\varphi_0(\ye,\ze)=O(\ke)\hbox{ and }\ye=O(\ke)
\eequa
when $\eps\to 0$. Noting that
$\varphi_0(\ye,\ze)=\varphi_0(0,\ze)+O(|\ye|)$ when $\eps\to 0$, we get
that $\varphi_0(0,\ze)=O(\ke)$. These last equalities,
(\ref{lim:prop32:1}), (\ref{lim:prop32:2}) and \eqref{eq:lundi} prove
(\ref{claim:1}).\hfill$\Box$

\medskip\noindent We let
\bequa\label{def:lae}
\lae=-\frac{\xune}{\ke}>0,\;\theta_\eps=\frac{\ye}{\ke}\in\P \hbox{ and
}\re=-\frac{\varphi_0(0,\ze)}{\ke}.
\eequa
It follows from  (\ref{claim:1}) and \eqref{lim:prop32:2} that there 
exist $\lambda_0\geq 0$,
$\rho_0\in\rr$ and $\theta_0\in\P $ such that
\bequa\label{lim:lae}
\lim_{\eps\to 0}\lae=\lambda_0,\; \lim_{\eps\to
0}\theta_\eps=\theta_0\hbox{ and }\lim_{\eps\to 0}\re=\rho_0.
\eequa
We claim that $\re\geq 0$ for all $\eps>0$. Indeed, since $\Porth 
\cap\Omega=\emptyset$, there exists $\delta>0$ such that for all
$z\in\hbox{span}\{\vec{e}_{k+1},...,\vec{e}_n\}\cap B_{\delta}(0)$
\bequa\label{ineq:varphi0}
\varphi_0(0,z)\leq 0.
\eequa
The definition (\ref{def:lae}) of $\re$ yields that $\re\geq
0$ for all $\eps>0$. Note that it follows from (\ref{ineq:varphi0}) 
that there exists $C>0$ such that
\bequa\label{ineq:d:pi}
d(x,\partial\Omega)\leq C|\pi(x)|
\eequa
for all $x\in\Omega$.

\medskip\noindent{\bf Step \ref{sec:part1}.4:} From now on,  we let
$\bze=(0,0,\ze)$
for all $\eps>0$ where $\ze$ is defined in (\ref{def:xune:ye:ze}), and for
any $x\in \frac{U-\bze}{\ke}\cap\{x_1\leq 0\}$, we set
\bequa\label{def:ve:sec3}
\ve(x):=\frac{\ue\circ\varphi(\bze+\ke x)}{\ue(\xe)},
\eequa
where $\varphi$ is defined in (\ref{def:vphi:3}). It follows from
(\ref{def:lae}) that
\bequa\label{eq:v:nonzero}
\ve(-\lae,\theta_\eps,0)=1.
\eequa
As easily checked, for any $\eta\in C_c^\infty(\rn)$, we have that
$\eta\ve\in\hunrnm$
for all $\eps>0$.

\medskip\noindent{\it Step \ref{sec:part1}.4.1:} There exists
$v\in\hunrnm$ such that for any $\eta\in C_c^\infty(\rn)$,
$$\eta\ve\rightharpoonup \eta v$$
weakly in $\hunrnm$ when $\eps\to 0$. The proof is rather similar to what
was done in \cite{gr1} to which we refer for details.\\

\medskip\noindent{\it Step \ref{sec:part1}.4.2:} We claim that
$\lim_{\eps\to 0}\ve=v$
in $C^1_{loc}(\overline{\rnm})$, where $v\not\equiv 0$.\\

  Indeed, let $R>0$ and for any $i,j=1,...,n$, we let
$(\tge)_{ij}=(\partial_i\varphi(\bze+\ke x),\partial_j\varphi(\bze+\ke
x))$, where $(\cdot,\cdot)$ denotes the Euclidean scalar product on $\rn$.
We consider $\tge$ as a metric on $\rn$. We let
$$\Delta_{\tge}=-{\tge}^{ij}\left(\partial_{ij}-\Gamma_{ij}^k(\tge)\partial_k\right),$$
where $\tge^{ij}:=(\tge^{-1})_{ij}$ are the coordinates of the inverse of
the tensor $\tge$ and the $\Gamma_{ij}^k(\tge)$'s are the Christoffel
symbols of the metric $\tge$. With a change of variable and the definition
(\ref{def:ve:sec3}), equation (\ref{syst:ue}) rewrites as
\bequa\label{eq:ve:sec3}
\Delta_{\tge}\ve+\ke^2\ae\circ\varphi(\bze+\ke
x)\ve=\frac{\ve^{\crit-1-\pe}}{\left|\frac{\pi(\varphi(\bze+\ke
x))}{\ke}\right|^s}\hbox{ in }{\mathcal D}'(\{x_1<0\})
\eequa
for all $\eps>0$. It follows from the definition (\ref{def:vphi:3}) of
$\varphi$ and (\ref{ineq:varphi0}) that there exists $C_R>0$ such that
$|\pi(\varphi(\bze+\ke x))|\geq C_R \ke|\pi(x)|$
for all $x\in\rnm\cap B_R(0)$. With (\ref{def:me:xe}) and
(\ref{def:ve:sec3}), we get that $0\leq\ve\leq 1$. With the method used in
the Appendix, we get that $(\ve)_{\eps>0}$ converges in $C^1_{loc}(\rnm)$.
Since $\ve\rightharpoonup v$ weakly in $\hunrnm$ when $k\to +\infty$, we
get that
$\lim_{\eps\to 0}\ve=v\hbox{ in }C^1_{loc}(\overline{\rnm}).$
With (\ref{eq:v:nonzero}) and \eqref{lim:lae}, we get that 
$v(-\lambda_0,\theta_0,0)=1$, and in particular, $v\not\equiv 0$ and 
$\lambda_0>0$.\hfill$\Box$

\medskip\noindent{\it Step \ref{sec:part1}.4.3:} We claim that
$\Delta v=\frac{v^{\crit-1}}{|\pi(x)|^s}\hbox{ in }{\mathcal D}'(\rnm)$
and that
$$\int_{\rnm}|\nabla v|^2\,
dx=\mu_{s,\P }(\Omega)^{\frac{\crit}{\crit-2}}=\mu_{s,\P }(\rnm)^{\frac{\crit}{\crit-2}}.$$
Indeed, by passing to the weak limit $\eps\to 0$ in (\ref{eq:ve:sec3}), we
get that
$$\Delta v=\frac{v^{\crit-1}}{|\pi(x)-(\rho_0,0,0)|^s}\hbox{ in }{\mathcal
D}'(\rnm).$$
Testing this equality with $v\in \hunrnm$ and using the optimal
Hardy-Sobolev inequality (\ref{def:mus}), we get that
\beqn
\left(\int_{\rnm}|\nabla v|^2\,
dx\right)^{\frac{\crit-2}{\crit}}&=&\frac{\int_{\rnm}|\nabla v|^2\,
dx}{\left(\int_{\rnm}\frac{v^{\crit}}{|\pi(x)-(\rho_0,0,0)|^s}\,
dx\right)^{\frac{2}{\crit}}}\nonumber\\
&\geq& \frac{\int_{\rnm}|\nabla v|^2\,
dx}{\left(\int_{\rnm}\frac{v^{\crit}}{|\pi(x)|^s}\,
dx\right)^{\frac{2}{\crit}}}\geq\mu_{s,\P }(\rnm).\label{ineq:sec3:rho}
\eeqn
Here, we have used that $|\pi(x)-(\rho_0,0,0)|\geq |\pi(x)|$ since
$\rho_0\geq 0$ and $x_1< 0$ for all $x\in\rnm$. We then obtain that
\bequa\label{ineq:nrj:v:1}
\int_{\rnm}|\nabla v|^2\, dx\geq \mu_{s,\P }(\rnm)^{\frac{\crit}{\crit-2}}.
\eequa
Moreover, see for instance \cite{gr1}, we have that $\int_{\rnm}|\nabla
v|^2\, dx\leq \mu_{s,\P }(\rnm)^{\frac{\crit}{\crit-2}}$. We then get that
$$\int_{\rnm}|\nabla v|^2\,
dx=\mu_{s,\P }(\Omega)^{\frac{\crit}{\crit-2}}=\mu_{s,\P }(\rnm)^{\frac{\crit}{\crit-2}},$$
and that
\bequa\label{lim:re}
\lim_{\eps\to 0}\re=\rho_0=0\hbox{ and }\lim_{\eps\to 0}\me^{\pe}=1.
\eequa
For this last assertion, we refer to \cite{gr1}.\hfill$\Box$

\medskip\noindent Proposition \ref{prop:sec3:3} now follows from Steps
\ref{sec:part1}.1 to \ref{sec:part1}.4.\hfill$\Box$

We shall also need the following two claims for the next section

\medskip\noindent{\bf Step \ref{sec:part1}.5:} Under the hypothesis of
Proposition \ref{prop:sec3:3}, we have that
\bequa\label{eq:limReps}
\lim_{R\to +\infty}\lim_{\eps\to 0}\int_{\Omega\setminus
B_{R\ke}(\varphi(\bze))}\frac{\ue^{\crit-\pe}}{|\pi(x)|^s}\, dx=0.
\eequa
We omit the proof which is quite similar to \cite{gr1}.

\medskip\noindent{\bf Step \ref{sec:part1}.6:} We also claim that
\bequa\label{cv:out:0}
\lim_{\eps\to 0}\ue=0\hbox{ in
}C^1_{loc}(\overline{\Omega}\setminus\{x_0\}).
\eequa
   Indeed, for $\delta>0$, it follows from (\ref{eq:limReps}) that
$$\lim_{\eps\to 0}\int_{\Omega\setminus
B_\delta(x_0)}\frac{\ue^{\crit-1-\pe}(x)}{|\pi(x)|^s}\, dx=0.$$
Using the techniques in the Appendix of \cite{gr1,gr2}, we get that
$$\lim_{\eps\to 0}\Vert\ue\Vert_{L^p(\Omega\setminus B_\delta(x_0))}=0$$
for all $p\geq 1$, and
the method developed in this paper's Appendix, we get
(\ref{cv:out:0}).\hfill$\Box$

\section{Blow-up analysis, Part II}\label{sec:part2}
This section is devoted to the proof of the following strong pointwise
estimate.
\begin{prop}\label{prop:fund:est}
Let $\Omega$ be a smooth bounded domain of $\rn$, $n\geq 3$ and let $\P$ be a linear vector subspace of $\rn$ such that  $2\leq
\hbox{dim}_{\rr}\P \leq n-1$. Let $s\in (0,2)$ and assume that
(\ref{hyp:P:1})  holds. For $(\pe)_{\eps>0}$ in $[0,\crit-2)$ such 
that $\lim_{\eps\to 0}\pe=0$ and
$(\ae)_{\eps>0}$ as in (\ref{hyp:ae}), we consider $(\ue)_{\eps>0}\in
\huno\cap C^2(\overline{\Omega}\setminus \Porth  )$ such that (\ref{syst:ue}),
(\ref{hyp:nrj:min}) and (\ref{hyp:blowup}) hold. We let $x_0$, $\varphi$,
$(\me)_{\eps>0}$ and $(\bze)_{\eps>0}$ as in Proposition
\ref{prop:sec3:3}. Then, there exists $C>0$ such that
\bequa\label{eq:fund:est}
\ue(x)\leq C 
+C\frac{\me^{\frac{n}{2}}d(x,\partial\Omega)}{\left(\me^2+|x-\varphi(\bze)|^2\right)^{\frac{n}{2}}}
\eequa
and
\bequa\label{eq:fund:est:2}
|\nabla\ue(x)|\leq C
d(x,\partial\Omega)+C\frac{\me^{\frac{n}{2}}}{\left(\me^2+|x-\varphi(\bze)|^2\right)^{\frac{n}{2}}}
\eequa
for all $\eps>0$ and all $x\in\Omega$.
\end{prop}
\smallskip\noindent{\it Proof:} We take inspiration in \cite{dhr}. We proceed in five steps.

\medskip\noindent{\bf Step \ref{sec:part2}.1:} We claim that there exists
$C>0$ such
that
\bequa\label{ineq:est:1}
|\pi(x)|^{\frac{n-2}{2}}\ue(x)^{1-\frac{\pe}{\crit-2}}\leq C
\eequa
for all $\eps>0$ and all $x\in\Omega$.

Indeed if not, we let $\ye\in\Omega$ such that
\bequa\label{hyp:wpe}
|\pi(\ye)|^{\frac{n-2}{2}}\ue(\ye)^{1-\frac{\pe}{\crit-2}}=\sup_{x\in\Omega}|\pi(x)|^{\frac{n-2}{2}}\ue(x)^{1-\frac{\pe}{\crit-2}}\to
+\infty
\eequa
as $\eps\to 0$. We then let
\bequa\label{def:nueps}
\nu_\eps:=\ue(\ye)^{-\frac{2}{n-2}}\hbox{ and
}\ell_\eps:=\nu_\eps^{1-\frac{\pe}{\crit-2}}
\eequa
for all $\eps>0$. It follows from (\ref{hyp:wpe}) and (\ref{def:nueps})
that
\bequa\label{lim:ye:nueps}
\lim_{\eps\to 0}\nu_\eps=0\hbox{ and }\lim_{\eps\to
0}\frac{|\pi(\ye)|}{\ell_\eps}=+\infty,
\eequa
and  from (\ref{def:me:xe}) and (\ref{lim:re}) that
\bequa\label{lim:nu:pe}
\lim_{\eps\to 0}\nu_\eps^{\pe}=1.
\eequa
We also let
\bequa\label{def:gamma:eps}
\gamma_{\eps}:=|\pi(y_{\eps})|^{\frac{s}{2}}|u_{\eps}(y_{\eps})|^{\frac{2-\crit+\pe}{2}},
\eequa
for all $\eps>0$. It follows from (\ref{lim:ye:nueps}) that
\bequa\label{lim:bis}
\lim_{\eps\to 0}\frac{\gamma_\eps}{|\pi(\ye)|}=0.
\eequa

\smallskip\noindent{\it Case \ref{sec:part2}.1.1:} We assume first that,
up to a subsequence,
there exists $\rho>0$ such that
\bequa\label{hyp:case1:wpe}
\frac{d(\ye,\partial\Omega)}{\gamma_\eps}\geq 3\rho
\eequa
for all $\eps>0$. For any $x\in B_{2\rho}(0)$ and any $\eps>0$, we let
\bequa\label{def:we:wpe}
\we(x):=\nu_\eps^{\frac{n-2}{2}}\ue(\ye+\gamma_\eps x).
\eequa
Note that $\we$ is well defined thanks to (\ref{hyp:case1:wpe}). With
(\ref{hyp:wpe}) and (\ref{def:gamma:eps}), we get that
$$\left|\frac{\pi(\ye)}{|\pi(\ye)|}+\frac{\gamma_\eps}{|\pi(\ye)|}\pi(x)\right|^{\frac{n-2}{2}}\we(x)^{1-\frac{\pe}{\crit-2}}\leq
1.$$
In particular, with (\ref{lim:bis}), there exists $C_0>0$ such that
\bequa\label{bnd:we:wpe}
0\leq \we(x)\leq C_0
\eequa
for all $x\in B_{2\rho}(0)$ and all $\eps>0$. With (\ref{syst:ue}), we get
that
$$\Delta\we+\gamma_\eps^2 a_\eps(\ye+\gamma_\eps
x)\we=\frac{\we^{\crit-1-\pe}}{\left|\frac{\pi(\ye)}{|\pi(\ye)|}+\frac{\gamma_\eps}{|\pi(\ye)|}\pi(x)\right|^s}$$
for all $x\in B_{2\rho}(0)$ and all $\eps>0$. Since (\ref{lim:ye:nueps})
and (\ref{bnd:we:wpe}) hold, it follows from standard elliptic theory that
there exists $w\in C^1(B_{2\rho}(0))$ such that $w(0)=1$ and
\bequa\label{lim:we:wpe}
\lim_{\eps\to 0}\we=w
\eequa
in $C_{loc}^1(B_{2\rho}(0))$. Mimicking what was done in Step
\ref{sec:part1}.1, we get a contradiction.

\smallskip\noindent{\it Case \ref{sec:part2}.1.2:} We assume that
\bequa\label{hyp:case2:wpe}
\lim_{\eps\to 0}\frac{d(\ye,\partial\Omega)}{\gamma_\eps}=0.
\eequa
As in Step \ref{sec:part1}.1, we get a contradiction. We refer to
\cite{gr1} for proof in a similar context.

\medskip\noindent In both cases, we have contradicted (\ref{hyp:wpe}).
This proves (\ref{ineq:est:1}).\hfill$\Box$

\medskip\noindent{\bf Step \ref{sec:part2}.2:} This step is a slight
improvement of
(\ref{ineq:est:1}). We claim that
\bequa\label{ineq:est:2}
\lim_{R\to +\infty}\lim_{\eps\to 0}\sup_{x\in \Omega\setminus
B_{R\ke}(\varphi(\bze))}|\pi(x)|^{\frac{n-2}{2}}\ue(x)^{1-\frac{\pe}{\crit-2}}=0.
\eequa
The proof is similar to Step \ref{sec:part2}.1 and uses the techniques
developed in \cite{gr1}. We refer to Step \ref{sec:part2}.1 and \cite{gr1}
for the details.

\medskip\noindent{\bf Step \ref{sec:part2}.3:} We claim that for any
$\nu\in (0,1)$ and any $R>0$, there exists $C(\nu,R)>0$ such that
\bequa\label{estim:nu:1}
\ue(x)\leq
C(\nu,R)\cdot\left(\frac{\me^{\frac{n}{2}-\nu(n-1)}d(x,\partial\Omega)^{1-\nu}}{(\me^2+|x-\varphi(\bze)|^2)^{\frac{n(1-\nu)}{2}}}+d(x,\partial\Omega)^{1-\nu}\right)
\eequa
for all $x\in\Omega$ and all $\eps>0$.\par

%\smallskip\noindent{\it Proof of the Claim:} Indeed, since $\Delta$ is
%coercive on $\Omega$, we
Indeed, let $G$ be the Green's function for $\Delta$ in $\Omega$ with Dirichlet
boundary condition, and set $H_\eps(x)=-\partial_{\vec{n}} G(x,\varphi(\bze))$
for all $x\in\overline{\Omega}\setminus\{\varphi(\bze)\}$, where here 
$\vec{n}$ denotes
the
outward normal vector at $\partial\Omega$. It follows from Theorem 9.2 of
\cite{gr2} that $H_\eps\in
C^2(\overline{\Omega}\setminus\{\varphi(\bze)\})$, that
\bequa\label{eq:Ge}
\Delta H_\eps=0
\eequa
in $\Omega$ and that there exist $\delta_1,C_1>0$ such that

\bequa\label{ineq:Ge:1}
\frac{d(x,\partial\Omega)}{C_1|x-\varphi(\bze)|^{n}}\leq H_\eps(x)\leq
\frac{C_1
d(x,\partial\Omega)}{|x-\varphi(\bze)|^{n}}
\eequa
and --using (\ref{ineq:d:pi})-- that
\bequa\label{ineq:Ge:2}
\frac{|\nabla H_\eps(x)|}{H_\eps(x)}\geq
\frac{1}{C_1'd(x,\partial\Omega)}\geq
\frac{1}{C_1 |\pi(x)|}
\eequa
for all $x\in \Omega\cap B_{2\delta_1}(0)$. Let $\lambda_{1}>0$ be the
first eigenvalue of $\Delta$ on $\Omega$, and let
$\psi\in C^2(\overline{\Omega})$ be ``the first eigenfunction" in such a
way that
$$\left\{\begin{array}{ll}
\Delta\psi=\lambda_{1}\psi & \hbox{ in
}\Omega\\
\psi>0 & \hbox{ in }\Omega\\
\psi=0 & \hbox{ on }\partial\Omega.
\end{array}\right.$$
It follows from standard elliptic theory, Hopf's maximum principle 
and again \eqref{ineq:d:pi} that
there exists $C_2,\delta_2>0$ such that
\bequa\label{ppty:phi}
\frac{1}{C_2}d(x,\partial\Omega)\leq \psi(x)\leq C_2
d(x,\partial\Omega)\hbox{ and
}\frac{|\nabla\psi(x)|}{\psi(x)}\geq\frac{1}{C_2'
d(x,\partial\Omega)}\geq\frac{1}{C_2 |\pi(x)|}
\eequa
for all $x\in \Omega\cap B_{2\delta_2}(\varphi(\bze))$. We now consider
the operator
$$L_\eps=\Delta
+\left(\ae-\frac{\ue^{\crit-2-\pe}}{|\pi(x)|^s}\right).$$

\medskip\noindent{\it Step \ref{sec:part2}.3.1:} We claim that there exist
$\delta_0>0$ and $R_0>0$ such that for any $\nu\in (0,1)$ and any $R>R_0$,
$\delta\in (0,\delta_0)$, we have that

\bequa\label{ineq:LG}
L_\eps H_\eps^{1-\nu}>0,\hbox{ and }L_\eps\psi^{1-\nu}>0
\eequa
for all $x\in \Omega\cap B_\delta(\varphi(\bze))\setminus
\overline{B}_{R\ke}(\varphi(\bze))$ and
for all $\eps>0$ sufficiently small.

\noindent Indeed, with (\ref{eq:Ge}), we get that

\bequa\label{eq:Ge:2}
\frac{L_\eps H_\eps^{1-\nu}}{H_\eps^{1-\nu}}(x)=\ae(x)+\nu(1-\nu)\frac{|\nabla
H_\eps|^2}{H_\eps^2}(x)-\frac{\ue(x)^{\crit-2-\pe}}{|\pi(x)|^s}
\eequa
for all $x\in \Omega$ and all $\eps>0$. We let $\alpha>0$. It follows from
(\ref{ineq:est:2}) that there exists $R_0>0$ such that for any $R>R_0$, we
have that
$$|\pi(x)|^{2-s}|\ue(x)|^{\crit-2-\pe}< \alpha$$
for all $x\in (B_\delta(\varphi(\bze))\setminus
\overline{B}_{R\ke}(\varphi(\bze)))\cap\Omega$
and all $\eps>0$ small enough. With \eqref{hyp:ae}, (\ref{eq:Ge:2}) 
and \eqref{ineq:Ge:2}, we get that
for $\alpha>0$ and $\delta>0$ small enough, we have that
$$\frac{L_\eps H_\eps^{1-\nu}}{H_\eps^{1-\nu}}(x)>\frac{\nu(1-\nu)-\alpha
C_1^2-C_1^2|\pi(x)|^2|\ae(x)|}{C_1^2|\pi(x)|^2}>0$$
for all $x\in (B_\delta(\varphi(\bze))\setminus
\overline{B}_{R\ke}(\varphi(\bze)))\cap\Omega$
and all $\eps>0$ small enough. The proof of the second inequality of
(\ref{ineq:LG}) goes the same way.

\medskip\noindent{\it Step \ref{sec:part2}.3.2:} It follows from
(\ref{cv:v:c1}) in Proposition \ref{prop:sec3:3} that there exists
$C_1(R)>0$ such that
$$\ue(x)\leq C_1(R) \me^{-\frac{n}{2}}d(x,\partial\Omega)$$
for all $x\in\Omega\cap\partial B_{R\ke}(\varphi(\bze))$ and all $\eps>0$.
In particular, there exists $C(R)>0$ such that
\bequa\label{sup:ue:boundary:2}
\ue(x)\leq C(R) \me^{\frac{n}{2}-\nu(n-1)}H_\eps^{1-\nu}(x)
\eequa
for all $x\in \Omega\cap \partial B_{R\ke}(\varphi(\bze))$ and all
$\eps>0$.

\smallskip\noindent It follows from (\ref{cv:out:0}) there exists
$C_1(\delta)>0$ such that
\bequa\label{sup:ue:boundary:1}
\ue(x)\leq C_1(\delta)d(x,\partial\Omega)
\eequa
for all $x\in\Omega\cap\partial B_{\delta}(\varphi(\bze))$ and all
$\eps>0$. In particular, there exists $C(\delta)>0$ such that
$\ue(x)\leq C(\delta) \psi(x)^{1-\nu}$
for all $x\in \Omega\cap\partial B_\delta(\varphi(\bze))$ and all
$\eps>0$. We let
$$D_{\eps,R,\delta}:=(B_\delta(\varphi(\bze))\setminus
\overline{B}_{R\ke}(\varphi(\bze)))\cap\Omega.$$
It follows from (\ref{sup:ue:boundary:2}) and (\ref{sup:ue:boundary:1})
that
\bequa\label{control:ue:bord}
\ue(x)\leq C(R) \me^{\frac{n}{2}-\nu(n-1)}H_\eps^{1-\nu}(x)+C(\delta)
\psi(x)^{1-\nu}
\eequa
for all $\eps>0$ and all $x\in \partial D_{\eps,R,\delta}$.

\medskip\noindent{\it Step \ref{sec:part2}.3.3:} We claim that $L_\eps$ is
coercive and therefore verifies the comparison principle on
$D_{\eps,R,\delta}$. \\

   Indeed, with (\ref{eq:limReps}), we get that for any $\alpha>0$, there
exists $\tilde{R}_0>0$ such that for any $R>\tilde{R}_0$, we have that
$$\int_{\Omega\setminus
B_{R\ke}(\varphi(\bze))}\frac{\ue^{\crit-\pe}(x)}{|\pi(x)|^s}\, dx\leq
\alpha.$$
Since $\Delta+\ae$ is uniformly coercive, we get that $L_\eps$ is coercive
on $\Omega\setminus B_{R\ke}(\varphi(\bze))$ for $R$ large enough. We
refer to Lemma 3.4 of \cite{r4} for details on this assertion.

\medskip\noindent{\it Step \ref{sec:part2}.3.4:} Since
$$L_\eps(C(R) \me^{\frac{n}{2}-\nu(n-1)}H_\eps^{1-\nu}(x)+C(\delta)
\psi(x)^{1-\nu})>0=L_\eps\ue$$
in $D_{\eps,R,\delta}$ and (\ref{control:ue:bord}) holds, we get from Step
\ref{sec:part2}.3.3 that
$$\ue(x)\leq C(R) \me^{\frac{n}{2}-\nu(n-1)}H_\eps^{1-\nu}(x)+C(\delta)
\psi(x)^{1-\nu}$$
for all $x\in D_{\eps,R,\delta}$. With (\ref{ineq:Ge:1}) and
(\ref{ppty:phi}), we then get that
(\ref{estim:nu:1}) holds on
$D_{\eps,R,\delta}=(B_\delta(\varphi(\bze))\setminus
\overline{B}_{R k_{\eps}}(\varphi(\bze)))\cap\Omega$ for $R$ large and
$\delta$
small. It follows from this last assertion, (\ref{cv:v:c1}) in 
Proposition \ref{prop:sec3:3} and (\ref{cv:out:0}) that
(\ref{estim:nu:1}) holds on $\Omega$.\hfill$\Box$

\medskip\noindent{\bf Step \ref{sec:part2}.4:} We claim that there exists
$C>0$ such that
\bequa\label{ineq:ef}
\ue(x)\leq C
d(x,\partial\Omega)+C\frac{\me^{\frac{n}{2}}d(x,\partial\Omega)}{\left(\me^2+|x-\varphi(\bze)|^2\right)^{\frac{n}{2}}}
\eequa
for all $x\in\Omega$ and all $\eps>0$. \\

Indeed, it follows from (\ref{cv:v:c1}) in Proposition \ref{prop:sec3:3}
and (\ref{cv:out:0}) that for any $\delta, R>0$, inequality
(\ref{ineq:ef}) holds for all $x\in (\Omega\setminus
B_\delta(\varphi(\bze)))\cup (\Omega\cap B_{R\me}(\varphi(\bze)))$ for all
$\eps>0$. What is left is to prove  (\ref{ineq:ef}) for any sequence 
$(\ye)_{\eps>0}\in\Omega$ such
that
\bequa\label{hyp:ye}
\lim_{\eps\to 0}\ye=x_0\hbox{ and }\lim_{\eps\to
0}\frac{|\ye-\varphi(\bze)|}{\ke}=+\infty.
\eequa
We show that (\ref{ineq:ef}) holds for $x=\ye$. With Green's
representation formula, we get that
$$\ue(\ye)=\int_\Omega
G_\eps(\ye,y)\frac{\ue(y)^{\crit-1-\pe}}{|\pi(y)|^s}\, dy,$$
where $G_\eps$ is the Green's function for the uniformly coercive 
operator $\Delta+\ae$. For
$\nu\in (0,1)$, we use (\ref{estim:nu:1}) and (\ref{ineq:d:pi}) to get that
\beqn
\ue(\ye) &\leq & C\int_\Omega
G_\eps(\ye,y)\frac{d(y,\partial\Omega)^{(1-\nu)(\crit-1-\pe)}}{|\pi(y)|^s}\,
dy\nonumber\\
&&+C\int_\Omega\frac{G_\eps(\ye,y)}{|\pi(y)|^s}\left(\frac{\me^{\frac{n}{2}-(n-1)\nu} d(y,\partial\Omega)^{1-\nu}}{\left(\me^2+|y-\varphi(\bze)|^2\right)^{\frac{n(1-\nu)}{2}}}\right)^{\crit-1-\pe}\,dy\nonumber\\
&\leq & I_{\eps,1}+I_{\eps,2}+I_{\eps,3}\label{ineq:def:i}
\eeqn
where
$$I_{\eps,1}=\int_\Omega G_\eps(\ye,y)
|\pi(y)|^{(1-\nu)(\crit-1-\pe)-s}
\, dy,$$
$$I_{\eps,2}=\int_{D_{\eps,2}}\frac{G_\eps(\ye,y)}{|\pi(y)|^s}\left(\frac{\me^{\frac{n}{2}-(n-1)\nu}d(y,\partial\Omega)^{1-\nu}}{\left(\me^2+|y-\varphi(\bze)|^2\right)^{\frac{n(1-\nu)}{2}}}\right)^{\crit-1-\pe}\, dy,$$
and
$$I_{\eps,3}=\int_{D_{\eps,3}}\frac{G_\eps(\ye,y)}{|\pi(y)|^s}\left(\frac{\me^{\frac{n}{2}-(n-1)\nu}d(y,\partial\Omega)^{1-\nu}}{\left(\me^2+|y-\varphi(\bze)|^2\right)^{\frac{n(1-\nu)}{2}}}\right)^{\crit-1-\pe}\,dy$$
for all $\eps>0$, where
$$D_{\eps,2}:=\left\{|\ye-y|>\frac{1}{2}|\ye-\varphi(\bze)|\right\}\hbox{ and }D_{\eps,3}:=\left\{|\ye-y|<\frac{1}{2}|\ye-\varphi(\bze)|\right\}$$

\medskip\noindent We first deal with $I_{\eps,1}$. The Green's function
verifies
$$G_\eps(\ye,y)\leq C\frac{d(\ye,\partial\Omega)}{|\ye-y|^{n-1}}$$
for all $y\in \Omega\setminus\{\ye\}$ and all $\eps>0$. We refer to
\cite{gr2} for the proof of this assertion. Since $s\in (0,2)$ and
$\varphi(\bze)\in\partial\Omega$, we then get that
\bequa\label{ineq:i1}
I_{\eps,1}\leq
Cd(\ye,\partial\Omega)\int_{\Omega}\frac{|\pi(y)|^{(1-\nu)(\crit-1-\pe)-s}}{|\ye-y|^{n-1}}\leq C d(\ye,\partial\Omega)
\eequa
for all $\eps>0$.

\medskip\noindent For $I_{\eps,2}$, we note that the Green's function
verifies
\bequa\label{ineq:green:annexe}
G_\eps(\ye,y)\leq
C\frac{d(\ye,\partial\Omega)d(y,\partial\Omega)}{|\ye-y|^{n}}
\eequa
for all $y\in \Omega\setminus\{\ye\}$ and all $\eps>0$. We again refer to
\cite{gr2} for the proof of this assertion. We then get with
(\ref{ineq:d:pi}) and a change of variables that
\beqn
&&I_{\eps,2}\nonumber\\
&&\leq
C\int_{D_{\eps,2}}\frac{d(\ye,\partial\Omega)
}{|\ye-y|^n}\frac{\me^{(\frac{n}{2}-(n-1)\nu)(\crit-1-\pe)}
d(y,\partial\Omega)^{(1-\nu)(\crit-1-\pe)+1-s}}{\left(\me^2+|x-\varphi(\bze)|^2\right)^{\frac{n(1-\nu)}{2}(\crit-1-\pe)}}\, dy\nonumber\\
&&\leq
C\frac{d(\ye,\partial\Omega)\me^{(\frac{n}{2}-(n-1)\nu)(\crit-1-\pe)}}{|\ye-\varphi(\bze)|^n}\int_\Omega\frac{|y-\varphi(\bze)|^{(1-\nu)(\crit-1-\pe)+1-s}}{\left(\me^2+|x-\varphi(\bze)|^2\right)^{\frac{n(1-\nu)}{2}(\crit-1-\pe)}}\, dy\nonumber\\
&&\leq
C\frac{d(\ye,\partial\Omega)\me^{\frac{n}{2}}}{|\ye-\varphi(\bze)|^n}\int_{\rn}\frac{|z|^{(1-\nu)(\crit-1-\pe)+1-s}}{\left(1+|z|^2\right)^{\frac{n(1-\nu)}{2}(\crit-1-\pe)}}\, dy \leq
C\frac{d(\ye,\partial\Omega)\me^{\frac{n}{2}}}{|\ye-\varphi(\bze)|^n}.\label{ineq:i2}
\eeqn

\medskip\noindent To deal with $I_{\eps,3}$, we first note that for any $y\in D_{\eps,3}$, we have that
\bequa
\frac{1}{2}|\ye-\varphi(\bze)|\leq |y-\varphi(\bze)|\leq
\frac{3}{2}|\ye-\varphi(\bze)|.
\eequa
With inequality (\ref{ineq:green:app}) (with $\theta=1$) on the 
Green's function, we
then get that
\beq
&&I_{\eps,3}\\
&&\leq
C\frac{d(\ye,\partial\Omega)\me^{(\frac{n}{2}-(n-1)\nu)(\crit-1-\pe)}}{|\ye-\varphi(\bze)|^{n(1-\nu)(\crit-1-\pe)}}\int_{D_{\eps,3}}\frac{dy}{|\ye-y|^{n-1}|\pi(y)|^{s-(1-\nu)(\crit-1-\pe)}}.
\eeq
We let
$$\theta_\eps=\frac{\ye-\varphi(\bze)}{|\ye-\varphi(\bze)|}+\frac{(\varphi_0(0,\bze),0,0)}{|\ye-\varphi(\bze)|}.$$
With (\ref{def:lae}), (\ref{lim:re}) and (\ref{hyp:ye}), we get that there
exists $\theta_0\in\rn$ such that $|\theta_0|=1$ and $\lim_{\eps\to
0}\theta_\eps=\theta_0$. With the change of variables
$y=\ye+|\ye-\varphi(\bze)|z$ and using (\ref{hyp:ye}), we get that
\beqn
&&I_{\eps,3}\leq\\
&&
C\frac{d(\ye,\partial\Omega)\me^{(\frac{n}{2}-(n-1)\nu)(\crit-1-\pe)}}{|\ye-\varphi(\bze)|^{(n-1)(1-\nu)(\crit-1-\pe)+s-1}}\int_{|z|<\frac{1}{2}}\frac{dz}{|z|^{n-1}|\pi(\theta_\eps+z)|^{s-(1-\nu)(\crit-1-\pe)}}\nonumber\\
&&\leq 
C\frac{d(\ye,\partial\Omega)\me^{(\frac{n}{2}-(n-1)\nu)(\crit-1-\pe)}}{|\ye-\varphi(\bze)|^{(n-1)(1-\nu)(\crit-1-\pe)+s-1}}\nonumber\\
&&=
o\left(\frac{d(\ye,\partial\Omega)\me^{\frac{n}{2}}}{|\ye-\varphi(\bze)|^n}\right)\label{ineq:i3}
\eeqn
when $\eps\to 0$. Plugging (\ref{ineq:i1}), (\ref{ineq:i2}) and
(\ref{ineq:i3}) in (\ref{ineq:def:i}) and using again (\ref{hyp:ye}), we
get that
$$\ue(\ye)\leq C
d(\ye,\partial\Omega)+C\frac{\me^{\frac{n}{2}}d(\ye,\partial\Omega)}{\left(\me^2+|\ye-\varphi(\bze)|^2\right)^{\frac{n}{2}}}$$
when $\eps\to 0$. This ends the proof of (\ref{ineq:ef}).

\medskip\noindent{\bf Step \ref{sec:part2}.5:} We claim that there exists
$C>0$ such that
\bequa\label{ineq:fe:2}
|\nabla\ue(x)|\leq
C+C\frac{\me^{\frac{n}{2}}}{\left(\me^2+|x-\varphi(\bze)|^2\right)^{\frac{n}{2}}}
\eequa
for all $x\in\Omega$.

To prove the claim, as in Step \ref{sec:part2}.4, we just need to consider $(\ye)_{\eps>0}\in \Omega$ as in \eqref{hyp:ye}. We use Green's
representation formula to write
$$\nabla\ue(\ye)=\int_\Omega\nabla_xG_\eps(\ye,y)\frac{\ue(y)^{\crit-1-\pe}}{|\pi(y)|^s}\, dy.$$
With (\ref{ineq:ef}), we get that
\bequa\label{ineq:def:j}
|\nabla\ue(\ye)|\leq J_{\eps,1}+J_{\eps,2}+J_{\eps,3},
\eequa
where
$$J_{\eps,1}=C\int_\Omega|\nabla_xG_\eps(\ye,y)|\frac{d(y,\partial\Omega)^{\crit-1-\pe}}{|\pi(y)|^s}\, dy,$$
$$J_{\eps,2}=C\int_{|\ye-y|>\frac{1}{2}|\ye-\varphi(\bze)|}|\nabla_xG_\eps(\ye,y)|\frac{\me^{\frac{n}{2}(\crit-1-\pe)}d(y,\partial\Omega)^{\crit-1-\pe}}{|\pi(y)|^s\left(\me^2+|y-\varphi(\bze)|^2\right)^{\frac{n}{2}(\crit-1-\pe)}}dy$$
and
$$J_{\eps,3}=C\int_{|\ye-y|<\frac{1}{2}|\ye-\varphi(\bze)|}|\nabla_xG_\eps(\ye,y)|\frac{\me^{\frac{n}{2}(\crit-1-\pe)}d(y,\partial\Omega)^{\crit-1-\pe}}{|\pi(y)|^s\left(\me^2+|y-\varphi(\bze)|^2\right)^{\frac{n}{2}(\crit-1-\pe)}}dy.$$

\noindent To estimate $J_{\eps,1}$, use that the Green's function
satisfies
\bequa\label{ineq:green:uneautre}
|\nabla_x G(\ye,y)|\leq \frac{C}{|\ye-y|^{n-1}}
\eequa
for all $y\in\Omega\setminus\{\ye\}$ and all $\eps>0$. We refer to
\cite{gr2} for the proof of this inequality. With (\ref{ineq:d:pi}), we
then get that
\bequa\label{ineq:j1}
J_{\eps,1}\leq
C\int_\Omega\frac{dy}{|\ye-y|^{n-1}|\pi(y)|^{s-(\crit-1-\pe)}}\leq C
\eequa

\medskip\noindent For $J_{\eps,2}$, we use that (see \cite{gr2})
$$|\nabla_x G(\ye,y)|\leq \frac{Cd(y,\partial\Omega)}{|\ye-y|^{n}}$$
for all $y\in\Omega\setminus\{\ye\}$ and all $\eps>0$. Plugging this
inequality in $J_{\eps,2}$ and performing computations similar to what was
done in the proof of (\ref{ineq:i2}), we get that
\bequa\label{ineq:j2}
J_{\eps,2}\leq
C\frac{\me^{\frac{n}{2}}}{\left(\me^2+|x-\varphi(\bze)|^2\right)^{\frac{n}{2}}}
\eequa

\medskip\noindent To  deal finally with $J_{\eps,3}$, we again use
estimate (\ref{ineq:green:uneautre}) on the Green's function combined with
the same techniques as in the proof of (\ref{ineq:i3}), to obtain
\bequa\label{ineq:j3}
J_{\eps,3}\leq
C\frac{\me^{\frac{n}{2}}}{\left(\me^2+|x-\varphi(\bze)|^2\right)^{\frac{n}{2}}}
\eequa

\medskip\noindent Plugging (\ref{ineq:j1}), (\ref{ineq:j2}) and
(\ref{ineq:j3}) in (\ref{ineq:def:j}), we get
(\ref{ineq:fe:2}) and Proposition \ref{prop:fund:est}.\hfill$\Box$

\section{Pohozaev identity and proof of Theorem \ref{th:intro}}\label{sec:poho}
  We first  prove the following

\begin{prop}\label{prop:poho}
Let $\Omega$ be a smooth bounded domain of $\rn$, $n\geq 3$ and let $\P$ be a linear vector subspace of $\rn$ such that $2\leq
\hbox{dim}_{\rr}\P\leq n-1$. Assume that $s\in (0,2)$ and  that
(\ref{hyp:P:1})  holds. For $(\pe)_{\eps>0} \in [0,\crit-2)$ and
$(\ae)_{\eps>0}$ as in (\ref{hyp:ae}), we consider $(\ue)_{\eps>0}\in
\huno\cap C^2(\overline{\Omega}\setminus \Porth  )$ such that (\ref{syst:ue}),
(\ref{hyp:nrj:min}) and (\ref{hyp:blowup}) hold. Then there exist
$x_0\in\partial\Omega\cap \Porth  $,  $\gamma_0\geq 0$ and  a family
$(\me)_{\eps>0}\in \rr_{+}$ such that $\lim_{\eps\to 0}\me=0$ and
\bequa\label{lim:pe:ke}
\lim_{\eps\to
0}\frac{\pe}{\me}=\frac{2(n-s)}{(n-2)^2}\mu_{s,\P }(\rnm)^{-\frac{n-s}{2-s}}\int_{\partial\rnm}\left(\frac{1}{2}II_{x_0}(x,x)-\gamma_0\right)|\nabla v|^2\, dx,
\eequa
where $II_{x_0}$ is the second fundamental form of $\partial \Omega$ at
$x_0$.
\end{prop}

Sections \ref{sec:poho}.1 to \ref{sec:poho}.3 below are devoted to the proof of Proposition \ref{prop:poho}, while Theorem \ref{th:intro} and Corollary \ref{th:eq} are proved in Step \ref{sec:poho}.4.

\medskip\noindent{\bf Step \ref{sec:poho}.1:} We establish a Pohozaev-type
identity for $\ue$. In the sequel, we let $(\bze)_{\eps>0}$,
$(\me)_{\eps>0}$, $(\ke)_{\eps>0}$ and $x_0\in \Porth  \cap \partial\Omega$ as
in Proposition \ref{prop:fund:est}. We also consider the chart $\varphi$
defined in (\ref{def:vphi:3}). We let
$$V_\eps=\Omega\cap \varphi(B_{\sqrt{\me}}(\bze))=\varphi(\rnm\cap
B_{\sqrt{\me}}(\bze)).$$
In particular,
$$\partial V_\eps=\varphi(\rnm\cap \partial B_{\sqrt{\me}}(\bze))\cup
\varphi(B_{\sqrt{\me}}(\bze)\cap\partial\rnm)= V_\eps^1\cup V_\eps^2.$$
In the sequel, we denote by $\nu(x)$ the outward normal vector at
$x\in\partial V_\eps$ of the oriented hypersurface $\partial V_\eps$ (this
is defined outside a null measure set). Let $\tilde{x}_0\in\rn$. After
integrations by parts (for instance, we refer to \cite{gr1,gr2}), we get that
\beqn
&&\left(\frac{n-2}{2}-\frac{n-s}{\crit-\pe}\right)\int_{V_\eps}\frac{\ue^{\crit-\pe}}{|\pi(x)|^s}\, dx-s\int_{V_\eps}\frac{(\tilde{x}_0,\pi(x))}{|\pi(x)|^{s+2}}\cdot\frac{\ue^{\crit-\pe}}{\crit-\pe}\, dx\nonumber\\
&&+\int_{V_\eps}\left(\ae+\frac{(x-\tilde{x}_0)^i\partial_i\ae}{2}\right)\ue^2\, 
dx\nonumber\\
&&=\int_{\partial
V_\eps}\left(-\frac{n-2}{2}\ue\partial_\nu\ue+(x-\tilde{x}_0,\nu)\frac{|\nabla\ue|^2}{2}-(x-\tilde{x}_0)^i\partial_i\ue\partial_{\nu}\ue\right.\nonumber\\
&&\left.-\frac{(x-\tilde{x}_0,\nu)}{\crit-\pe}\cdot\frac{\ue^{\crit-\pe}}{|\pi(x)|^s}+\frac{\ae(x-\tilde{x}_0,\nu)}{2}\ue^2\right)\, d\sigma\label{eq:poho:4:1}
\eeqn
for all $\eps>0$. Since $\ue\equiv 0$ on $\partial\Omega$, taking
$\tilde{x}_0=\varphi(\bze)$ in (\ref{eq:poho:4:1}), we get that
\beqn
&&\left(\frac{n-2}{2}-\frac{n-s}{\crit-\pe}\right)\int_{V_\eps}\frac{\ue^{\crit-\pe}}{|\pi(x)|^s}\, dx-s\int_{V_\eps}\frac{(\varphi(\bze),\pi(x))}{|\pi(x)|^{s+2}}\cdot\frac{\ue^{\crit-\pe}}{\crit-\pe}\, dx\nonumber\\
&&+\int_{V_\eps}\left(\ae+\frac{(x-\varphi(\bze))^i\partial_i\ae}{2}\right)\ue^2\,dx\nonumber\\
&&=\int_{V_\eps^1}\left(-\frac{n-2}{2}\ue\partial_\nu\ue+(x-\varphi(\bze),\nu)\frac{|\nabla\ue|^2}{2}\right.\nonumber\\
&&\left.-(x-\varphi(\bze))^i\partial_i\ue\partial_{\nu}\ue-\frac{(x-\varphi(\bze),\nu)}{\crit-\pe}\cdot\frac{\ue^{\crit-\pe}}{|\pi(x)|^s}+\frac{\ae(x-\varphi(\bze),\nu)}{2}\ue^2\right)\,
d\sigma\nonumber\\
&&-\frac{1}{2}\int_{V_\eps^2}(x-\varphi(\bze),\nu)|\nabla\ue|^2\,
d\sigma.
\label{eq:poho:4}
\eeqn
With (\ref{hyp:nrj:min}), (\ref{eq:limReps}), (\ref{eq:fund:est}), 
(\ref{eq:fund:est:2}) and Proposition \ref{prop:sec3:3}, we get that
\beqn
&&\left(\frac{(n-2)^2}{4(n-s)}\mu_{s,\P }(\rnm)^{\frac{\crit}{\crit-2}}+o(1)\right)\pe+s\int_{V_\eps}\frac{(\varphi(\bze),\pi(x))}{|\pi(x)|^{s+2}}\cdot\frac{\ue^{\crit-\pe}}{\crit-\pe}\,dx\nonumber\\
&&=\frac{1}{2}\int_{V_\eps^2}(x-\varphi(\bze),\nu)|\nabla\ue|^2\,
d\sigma+o(\me)\label{eq:poho:5}.
\eeqn

\medskip\noindent{\bf Step \ref{sec:poho}.2:} We deal with the RHS of
(\ref{eq:poho:5}). With a change of variable, we get that
\beqn
&&\int_{\varphi(B_{\sqrt{\me}}(\bze)\cap\partial\rnm)}(x-\varphi(\bze),\nu)|\nabla\ue|^2\,d\sigma=\label{eq:poho:rhs:1}\\
&&(1+o(1))\me \int_{D_\eps}\left(\frac{\varphi(\bze+\ke
x)-\varphi(\bze)}{\ke^2},\nu\circ\varphi(\bze+\ke
x)\right)|\nabla\ve|_{\tge}^2\sqrt{|\tge|}\, dx\nonumber
\eeqn
where the metric $\tge$ is such that
$(\tge)_{ij}=(\partial_i\varphi,\partial_j\varphi)(\bze+\ke x)$ for all
$i,j=2,...,n$, $\ve$ is as in Proposition \ref{prop:sec3:3} and
$$D_\eps=B_{\frac{\sqrt{\me}}{\ke}}(0)\cap\{x_1=0\}.$$
Using the expression of $\varphi$ (see (\ref{def:vphi:3})), we get (see
\cite{gr1,gr2} for details) that
\beqn
&&\left(\frac{\varphi(\bze+\ke
x)-\varphi(\bze)}{\ke^2},\nu\circ\varphi(\bze+\ke x)\right)\nonumber\\
&&=\frac{1+o(1)}{\ke^2}\left(\varphi_0(\bze+\ke
x)-\varphi_0(\bze)-\ke\sum_{i=2}^nx^i\partial_i\varphi(\bze+\ke
x)\right)\nonumber\\
&&=-\frac{1}{2}\sum_{i,j=2}^n\partial_{ij}\varphi(\bze)x^ix^j+o_\eps(1)|x|^2\label{dev:PS}
\eeqn
for $\eps>0$ and $x\in D_\eps$. In this expression, $\lim_{\eps\to 
0}o_\eps(1)=0$ uniformly in $D_\eps$. Plugging
(\ref{dev:PS}) into (\ref{eq:poho:rhs:1}), using the estimate
(\ref{eq:fund:est:2}), Lebesgue's convergence theorem and 
\eqref{cv:v:c1}, we get that
\bequa\label{eq:poho:rhs:2}
\lim_{\eps\to
0}\frac{1}{\me}\int_{\varphi(B_{\sqrt{\me}}(\bze)\cap\partial\rnm)}(x-\varphi(\bze),\nu)|\nabla\ue|^2\, d\sigma=-\frac{1}{2}\int_{\partial\rnm}\partial_{ij}\varphi_0(0)x^ix^j|\nabla
v|^2\, dx.
\eequa

\medskip\noindent{\bf Step \ref{sec:poho}.3:} We deal with the second term
of the LHS of (\ref{eq:poho:5}). With the pointwise estimate
(\ref{eq:fund:est}) and a change of variables, we get that
\beq
&&\int_{V_\eps}\frac{(\varphi(\bze),\pi(x))}{|\pi(x)|^{s+2}}\cdot\frac{\ue^{\crit-\pe}}{\crit-\pe}\, dx\\
&&=\frac{1+o(1)}{\me^2}\int_{D'_\eps}\frac{(\pi\circ\varphi(\bze),\pi\circ\varphi(\bze+\ke x))}{\left|\frac{\pi\circ\varphi(\bze+\ke
x)}{\ke}\right|^{s+2}}\cdot\frac{\ve^{\crit-\pe}}{\crit-\pe}\, dx
\eeq
when $\eps\to 0$, where
$$D'_\eps=B_{R\frac{\sqrt{\me}}{\ke}}(0)\cap\{x_1<0\}.$$
With the explicit expression of $\varphi$ (see (\ref{def:vphi:3})) 
and noting $x=(x_1,y,z)$ as in \eqref{def:vphi:3}, we get
that
\beq
&&\int_{V_\eps}\frac{(\varphi(\bze),\pi(x))}{|\pi(x)|^{s+2}}\cdot\frac{\ue^{\crit-\pe}}{\crit-\pe}\, dx\\
&&=(1+o(1))\frac{\varphi_0(0,\ze)}{\ke}\int_{D'_\eps}\frac{x_1+\frac{\varphi_0(\ke y,\ze+\ke z))}{\ke}}{\left|\pi\left(x_1+\frac{\varphi_0(\ke y,\ze+\ke
z)}{\ke},y,z\right)\right|^{s+2}}\frac{\ve^{\crit-\pe}}{\crit-\pe}\, dx
\eeq
when $\eps\to 0$. With point (iii) of (\ref{def:vphi:3}), the estimate
(\ref{eq:fund:est}) and Lebesgue's convergence theorem, we get that
\beqn
&&\int_{V_\eps}\frac{(\varphi(\bze),\pi(x))}{|\pi(x)|^{s+2}}\frac{\ue^{\crit-\pe}}{\crit-\pe}\, dx\nonumber\\
&&=\frac{\varphi_0(0,\ze)}{\me}\left(\int_{\rnm}\frac{x_1 v^{\crit}}{\crit
|\pi(x)|^{s+2}}\, dx+o(1)\right)\label{eq:poho:lhs:1}
\eeqn
where $\lim_{\eps\to 0}o(1)=0$. Plugging (\ref{eq:poho:rhs:2}) and
(\ref{eq:poho:lhs:1}) into (\ref{eq:poho:5}) and noting that
$\varphi_0(0,\ze)\leq 0$ (see \eqref{ineq:varphi0}), we get that
\beqn
&&\left(\frac{(n-2)^2}{4(n-s)}\mu_{s,\P }(\rnm)^{\frac{n-s}{2-s}}+o(1)\right)\pe+\left(\int_{\rnm}\frac{s|x_1| v^{\crit}}{\crit |\pi(x)|^{s+2}}\, dx+o(1)\right)
\frac{|\varphi_0(0,\ze)|}{\me}\nonumber\\
&&=\left(-\frac{1}{4}\int_{\partial\rnm}\partial_{ij}\varphi_0(0)x^ix^j|\nabla v|^2\, dx+o(1)\right)\cdot\me\label{eq:poho:7}
\eeqn
where $\lim_{\eps\to 0}o(1)=0$. In particular, we get that
$|\varphi_0(0,\ze)|=O(\me^2)$ when $\eps\to 0$. We let
$$\gamma_0=-\lim_{\eps\to 0}\frac{\varphi_0(0,\ze)}{\me^2}\geq 0.$$
With (\ref{eq:poho:7}), we get that
\beq
&&\lim_{\eps\to
0}\frac{(n-2)^2}{4(n-s)}\mu_{s,\P }(\rnm)^{\frac{n-s}{2-s}}\frac{\pe}{\me}\\
&&=-\frac{1}{4}\int_{\partial\rnm}\partial_{ij}\varphi_0(0)x^ix^j|\nabla
v|^2\, dx -\gamma_0\frac{s}{\crit}\int_{\rnm}\frac{|x_1|\cdot
v^{\crit}}{|\pi(x)|^{s+2}}\, dx.
\eeq
Taking $\tilde{x_0}=\vec{e}_1$ in (\ref{eq:poho:4:1}), using a change of
variable and the arguments used to prove (\ref{eq:poho:lhs:1}), we get
that
\bequa\label{ipp:1}
\frac{s}{\crit}\int_{\rnm}\frac{|x_1| v^{\crit}}{ |\pi(x)|^{s+2}}\,
dx=\frac{1}{2}\int_{\partial\rnm}|\nabla v|^2\, dx.
\eequa
We consider the second fondamental form associated to $\partial\Omega$,
namely
$$II_p(x,y)=(d\nu_px,y)$$
for all $p\in\partial\Omega$ and all $x,y\in T_{x_0}\partial\Omega$
(recall that $\nu$ is the outward normal vector at the hypersurface
$\partial\Omega$). In the basis $(\vec{e}_1,...,\vec{e}_n)$, the matrix of
the bilinear form $II_{x_0}$ is $-D^2_0\varphi_0$, where $D^2_0\varphi_0$
is the Hessian matrix of $\varphi_0$ at $0$. With this remark,
(\ref{eq:poho:7}) and (\ref{ipp:1}), we get that
$$\lim_{\eps\to
0}\frac{\pe}{\me}=\frac{2(n-s)}{(n-2)^2}\mu_{s,\P }(\rnm)^{-\frac{n-s}{2-s}}\int_{\partial\rnm}\left(\frac{1}{2}II_{x_0}(x,x)-\gamma_0\right)|\nabla v|^2\, dx,$$
where $\gamma_0\geq 0$. This ends the proof of Proposition \ref{prop:poho}.

\medskip\noindent{\bf Step \ref{sec:poho}.4:} We are now in position to prove Theorem \ref{th:intro}. Points (A) and (B) of Theorem \ref{th:intro} are direct consequences of Propositions \ref{prop:Case1} and \ref{prop:Case2}.

To establish Part (C) of Theorem \ref{th:intro}, assume that
(\ref{hyp:P:1}) holds and let us suppose that there are no extremals for (\ref{def:mus}). It follows from Proposition \ref{prop:subcrit} that there exists $(\ue)_{\eps>0}\in\huno$ such that (\ref{syst:ue}) and (\ref{hyp:nrj:min}) hold with $\pe=\eps$ and $\ae\equiv 0$. Since there are no extremals, it follows from Proposition \ref{prop:subcrit} that (\ref{hyp:blowup}) holds. We apply Proposition \ref{prop:poho} and we get that
$$\lim_{\eps\to
0}\frac{\eps}{\me}=\frac{2(n-s)}{(n-2)^2}\mu_{s,\P }(\rnm)^{-\frac{n-s}{2-s}}\int_{\partial\rnm}\left(\frac{1}{2}II_{x_0}(x,x)-\gamma_0\right)|\nabla v|^2\, dx$$
where $x_0\in \Porth \cap\partial\Omega$ and $\gamma_0\geq 0$. We then get that
\bequa\label{lim:eps:me}
\int_{\partial\rnm}II_{x_0}(x,x)|\nabla v|^2\, dx\geq 0
\eequa
Assume that we are in the first case of point (C) of Theorem \ref{th:intro}. We
then get that $II_{x_0}(x,x)\leq 0$ for all $x\in\partial\rnm$, but
$II_{x_0}(x,x)\not\equiv 0$. A contradiction with (\ref{lim:eps:me}).

\smallskip\noindent To relate our main result to conditions on the mean
curvature, we now assume that
$\P \cap T_x\partial\Omega$ and $\Porth  $ are orthogonal for the 
bilinear form $II_{x_0}$, we get in the coordinates (\ref{choice:basis}) and the chart (\ref{def:vphi:3}) that $(II_{x_0})_{ij}=0$ when $i\in \{2,...,k\}$ and $j\in \{k+1,n\}$. In
particular, we have with (\ref{lim:eps:me}) that
\bequa\label{lim:eps:me:2}
\left(\sum_{i,j=2}^k(II_{x_0})_{ij}\int_{\partial\rnm}x^ix^j|\nabla v|^2\,
dx\right)+\left(\sum_{i,j=k+1}^n(II_{x_0})_{ij}\int_{\partial\rnm}x^ix^j|\nabla
v|^2\, dx\right)\geq 0.
\eequa
%Since $\P \cap T_x\partial\Omega$ and $\Porth  $ are orthogonal for the 
%bilinear form $II_{x_0}$,
The matrix of the second fundamental form of $\partial\Omega\cap(\Porth + (T_{x_0}\partial\Omega)^{\bot})$ 
at $x_0$ with respect to a given vector $\vec{X}$ is 
$\left((II_{x_0}(\vec{X}))_{ij}\right)_{i,j\geq
k+1}=-\left(\partial_{ij}\varphi_0(0)X^1\right)_{i,j\geq k+1}$. Since
$\nabla\varphi_0(0)=0$ and $\varphi_0(0,z)\leq 0$ for $z$ close to $0$, we
get that for any direction $\vec{X}$, the principal curvatures of $\partial\Omega\cap(\Porth + (T_{x_0}\partial\Omega)^{\bot})$ 
at $x_0$ have a sign. If the mean curvature 
vector of $\partial\Omega\cap(\Porth + (T_{x_0}\partial\Omega)^{\bot})$ 
at $x_0$
is assumed to be null, it then follows that the second fundamental form of $\partial\Omega\cap(\Porth + (T_{x_0}\partial\Omega)^{\bot})$ at $x_0$ is null, and we get then from
(\ref{lim:eps:me:2}) that 
\bequa\label{lim:eps:me:3}
\sum_{i,j=2}^k(II_{x_0})_{ij}\int_{\partial\rnm}x^ix^j|\nabla v|^2\,
dx\geq 0.
\eequa
Here, $v\in\hunrnm$ is positive, verifies $\Delta
v=\frac{v^{\crit-1}}{|\pi(x)|^s}$ weakly and that $v(x)\leq
C(1+|x|^2)^{-n}$ for all $x\in\rnm$ (this last statement is a consequence of (\ref{cv:v:c1}) and (\ref{eq:fund:est})). It follows from Proposition \ref{prop:sym} that there exists $\tilde{v}$ such that $v(x_1,y,z)=\tilde{v}(x_1,|y|,z)$.
With this symmetry property, we get with (\ref{lim:eps:me:3}) that
$\sum_{i=1}^k(II_{x_0})_{ii}\geq 0$, and then the mean curvature at $x_0$ of $\partial\Omega$ is nonnegative. A contradiction with the assumption (2) of case (C) of Theorem \ref{th:intro}. This ends the proof of the Theorem.

\medskip\noindent Concerning Corollary \ref{th:eq}, the subcritical problem yields families of positive solutions to (\ref{syst:ue}) and (\ref{hyp:nrj:min}) with $\ae\equiv a$ and $\pe=\eps$. The proof of the result then goes as in the Proof of Theorem \ref{th:intro}.

\section{Proof of Proposition \ref{prop:n-1}}
We let $\Omega$ and $\P$ as in Proposition \ref{prop:n-1}. In particular $\dim_{\rr}\P =1$. The proof of case (B) of Proposition \ref{prop:n-1} goes exactly as the proof of Proposition \ref{prop:Case2}. Concerning case (C), we claim that $\mu_{s,\P }(\Omega)=0$ when $s\in [1,2)$. Indeed, taking $u\in C_c^\infty(\Omega)$ such that $u(x_0)=1$, where $x_0\in \Porth \cap\Omega$, it is easily checked that $\int_\Omega\frac{u^{\crit}}{|\pi(x)|^s}\, dx=+\infty$, and then
$\mu_{s,\P }(\Omega)=0$ is not achieved. When $s\in (0,1)$ in case (A),
%\footnote{Is inequality  \eqref{ineq:HS:rn}  still valid?},
the proof of non-achievement goes as the proof of Proposition \ref{prop:Case1}.

\medskip\noindent We are left with case (C) of Proposition \ref{prop:n-1},
that is $\Porth \cap\Omega=\emptyset$ and $\Porth \cap\partial\Omega\neq\emptyset$.
Up to a change of coordinates, we assume that $\Porth =\{x_1=0\}$,
$\Omega\subset\rnm$ and $0\in \Porth \cap\partial\Omega$ and $|\pi(x)=|x_1|$. In particular, it follows from the
Sobolev inequality and the Hardy inequality
% (see for instance \cite{ms}\footnote{reference to check})
that
$$\mu_{s,\P }(\rnm):=\inf\left\{\frac{\int_{\rnm}|\nabla u|^2\,
dx}{\left(\int_{\rnm}\frac{|u|^{\crit}}{|x_1|^s}\,
dx\right)^{\frac{2}{\crit}}}; u\in \hunrnm\setminus\{0\}\right\}>0.$$
Since $\Omega\subset\rnm$, we get that $\mu_{s,\P }(\Omega)\geq
\mu_{s,\P }(\rnm)$. With arguments similar to the proof of Proposition \ref{prop:test}, we also get the reverse inequality, and then
$\mu_{s,\P }(\Omega)=\mu_{s,\P }(\rnm)$. In particular, an extremal for $\mu_{s,\P }(\Omega)$ is an extremal for $\mu_{s,\P }(\rnm)$ and vice-versa. As in the proof of Proposition \ref{prop:Case1}, the maximum principle yields a contradiction.

\section{Appendix: Regularity of weak solutions}\label{sec:app}
In this appendix, we prove the following regularity result:
\begin{prop}\label{prop:app}
Let $\Omega$ be a smooth bounded domain of $\rn$, $n\geq 3$. Let
$\P\subset\rn$ be a $k-$dimensional linear subspace of $\rn$, where $2\leq k\leq n-1$. We assume that
$$\Porth  \cap\Omega=\emptyset\hbox{ and }\Porth  \cap\partial\Omega\neq\emptyset.$$
We let $s\in (0,2)$ and $a\in C^{0,\alpha}(\overline{\Omega})$, where
$\alpha\in (0,1)$. We let $\eps\in [0,\crit-2)$ and consider $u\in\huno$ a
weak solution of
\bequa\label{eq:u:app}
\Delta u+au=\frac{|u|^{\crit-2-\eps}u}{|\pi(x)|^s}  \hbox{ in }{\mathcal
D}'(\Omega).
\eequa
Then $u\in C^1(\overline{\Omega})\cap
C^{2,\alpha}(\overline{\Omega}\setminus\Porth)$.
\end{prop}

\medskip\noindent{\it Proof of Proposition \ref{prop:app}:}
Note that since $\crit<\frac{2n}{n-2}$, it follows from standard elliptic
theory that $u\in C^{2,\alpha}(\overline{\Omega}\setminus\Porth  )$. In
particular, $u\in C^{2,\alpha}(\Omega)$.

\medskip\noindent{\bf Step \ref{sec:app}.1:} We claim that
\bequa\label{bnd:lp}
u\in L^p(\Omega)
\eequa
for all $p\geq 1$. Indeed, the proof is similar to the case $\P=\rn$
provided in \cite{gr1,gr2}. We omit the proof and refer to \cite{gr1,gr2}
for the details.

\medskip\noindent In particular, we get that
$\frac{|u|^{\crit-2-\eps}u}{|\pi(x)|^s}\in L^p(\Omega)$ for all $1\leq
p<\frac{k}{s}$. In the case $k=n$, we take $p>\frac{n}{2}$, and then $u\in
L^\infty(\Omega)$. A bootstrap argument (see also \cite{eg}) then yields that $u\in
C^1(\overline{\Omega})$. However, in the general case $2\leq k\leq n-1$,
such an argument using standard elliptic theory does not hold, and we have
to use the Green's function to prove the proposition.

\medskip\noindent{\bf Step \ref{sec:app}.2:} We let $\theta\in
(0,\min\{2-s,1\})$. We claim that there exists $C>0$ such that
\bequa\label{app:claim1}
|u(x)|\leq C d(x,\partial\Omega)^\theta
\eequa
for a.e. $x\in\Omega$.

\smallskip\noindent{\it Proof of the claim:} We let $(\eta_k)_{k\in\nn}\in
C^\infty_c(\Omega)$ such that $0\leq \eta_k\leq 1$ for all $k$ and 
$\eta_k(x)=1$ for $d(x,\partial\Omega)\geq 2k^{-1}$. We let
$(u_k)_{k\in\nn}\in H_{1,0}^2(\Omega)$ such that
\bequa\label{eq:uk}
\Delta u_k=\eta_k\left(\frac{|u|^{\crit-2-\eps}u}{|\pi(x)|^s}-au\right).
\eequa
Since $u\in C^2(\Omega)$ and $\Omega\cap \Porth  =\emptyset$, we get that
$u_k\in C^2(\overline{\Omega})$ for all $k\in\nn$. We let $G$ be the
Green's function for $\Delta$ with Dirichlet boundary condition. It
follows from Green's representation formula that
\bequa\label{app:u:green}
u_k(x)=\int_\Omega
G(x,y)\eta_k(y)\left(\frac{|u|^{\crit-1-\eps}(y)}{|\pi(y)|^s}-au\right)\,
dy
\eequa
for a.e. $x\in\Omega$. It follows from Theorem 9.1 of \cite{gr2} that
there exists $C>0$ such that
\bequa\label{ineq:green:app}
0<G(x,y)\leq C\frac{d(x,\partial\Omega)^\theta}{|x-y|^{n-2+\theta}}
\eequa
for all $x,y\in\Omega$, $x\neq y$. Plugging this inequality in
(\ref{app:u:green}) and using H\"older's inequality, we get that
\beqn
|u_k(x)|&\leq & C
d(x,\partial\Omega)^\theta\int_\Omega\frac{1}{|x-y|^{n-2+\theta}}\left(\frac{|u(y)|^{\crit-1-\eps}}{|\pi(y)|^s}+|u(y)|\right)\, dy\nonumber\\
&\leq & C d(x,\partial\Omega)^\theta \Vert
|u|^{\crit-1-\eps}\Vert_q\left(\int_\Omega\frac{dy}{|x-y|^{p(n-2+\theta)}|\pi(y)|^{sp}}\right)^{\frac{1}{p}}\nonumber\\
&&+ C d(x,\partial\Omega)^\theta \Vert
u\Vert_{q'}\left(\int_\Omega\frac{dy}{|x-y|^{p'(n-2+\theta)}}\right)^{\frac{1}{p'}}\label{ineq:uk}
\eeqn
where $p,q,p',q'>1$ are such that
$\frac{1}{p}+\frac{1}{q}=\frac{1}{p'}+\frac{1}{q'}=1$. Since $\theta\in
(0,1)$ and (\ref{bnd:lp}) holds, we get that there exists $C>0$ such that
for $p,p'>1$ sufficiently close to $1$, we have that
$$|u(x)|\leq
Cd(x,\partial\Omega)^\theta\left(\int_\Omega\frac{dy}{|x-y|^{p(n-2+\theta)}|\pi(y)|^{sp}}\right)^{\frac{1}{p}}+ C d(x,\partial\Omega)^\theta$$
for all $x\in\Omega$. For simplicity, up to a change of coordinates, we
write any $y\in\rn$ as $y=(y',y'')$, where $y'=\pi(y)\in \rr^k=\P $ and
$y''\in \rr^{n-k}=\Porth $. We let $R>0$ such that $\Omega\subset B_R^k(0)\times
B_R^{n-k}(0)$ (the product of the ball of radius $R$ in $\rr^k$ and the
ball of radius $R$ in $\rr^{n-k}$). We then get with a change of variable that
\beq
&&\int_\Omega\frac{dy}{|x-y|^{p(n-2+\theta)}|\pi(y)|^{sp}}\\
&&\leq
C\int_{B_R^k(0)}\frac{1}{|y'|^{ps}}\int_{B_R(0)^{n-k}}\left(\frac{dy''}{|x'-y'|^{p(n-2+\theta)}+|x''-y''|^{p(n-2+\theta)}}\right)\, dy'\\
&&\leq
C\int_{B_R^k(0)}\frac{1}{|y'|^{ps}|x'-y'|^{p(n-2+\theta)+k-n}}\int_{B_{\frac{2R}{|x'-y'|}}(0)^{n-k}}\left(\frac{dz''}{1+|z''|^{p(n-2+\theta)}}\right)\, dy'\\
&&\leq
C\int_{B_R^k(0)}\frac{dy'}{|y'|^{ps}|x'-y'|^{p(n-2+\theta)+k-n}}\leq C
\eeq
for all $(x',x'')\in\Omega$. Here, we have taken $p>1$ close to $1$ and we
have used that $s\in (0,2)$. Plugging this inequality in (\ref{ineq:uk}),
we get that there exists $C>0$ such that
\bequa\label{bnd:uk}
|u_k(x)|\leq C d(x,\partial\Omega)^\theta
\eequa
for all $x\in\Omega$ and all $k\in\nn$. Multiplying (\ref{eq:uk}) by
$u_k$, integrating over $\Omega$, using that $u\in\huno$, the inequality
(\ref{ineq:HS:rn}) and (\ref{bnd:uk}), we get that there exists $C>0$ such
that $\Vert u_k\Vert_{\huno}\leq C$ for all $k\in\nn$. It then follows
that there exists $\tilde{u}\in\huno$ such that $u_k\rightharpoonup
\tilde{u}$ weakly in $\huno$ when $k\to +\infty$ and $\lim_{k\to
+\infty}u_k(x)=\tilde{u}(x)$ for a.e. $x\in\Omega$. The function
$\tilde{u}$ verifies $\Delta \tilde{u}=\frac{|u|^{\crit-2-\eps}}{|\pi(x)|^s}-au$ in ${\mathcal D}'(\Omega)$. Since $\Delta$ is coercive, it then follows from
(\ref{eq:u:app}) that $\tilde{u}=u$. With (\ref{bnd:uk}), we then get
(\ref{app:claim1}). \hfill$\Box$

\medskip\noindent{\bf Step \ref{sec:app}.3:} We claim that there exists
$C>0$ such that
\bequa\label{ineq:claim3}
|u(x)|\leq C d(x,\partial\Omega)
\eequa
for a.e. $x\in\Omega$.

\smallskip\noindent{\it Proof of the claim:} Indeed, we let $\theta_0\in
(0,1)$ such that there exists $C>0$ such that $|u(x)|\leq C
d(x,\partial\Omega)^{\theta_0}$. With (\ref{ineq:d:pi}), we get that there
exists $C>0$ such that $|u(x)|\leq C|\pi(x)|^{\theta_0}$ for all
$x\in\Omega$. We let $\theta\in (0,1)$. It follows from Green's
representation formula and (\ref{ineq:green:app}) that there exists $C>0$
such that
\beq
|u_k(x)|&=&\left|\int_\Omega
G(x,y)\eta_k(y)\left(\frac{|u|^{\crit-1-\eps}(y)}{|\pi(y)|^s}-au\right)\,
dy\right|\\
&\leq& C d(x,\partial\Omega)^\theta +C\int_\Omega
\frac{d(x,\partial\Omega)^{\theta}}{|x-y|^{n-2+\theta}|\pi(y)|^{s-\theta_0(\crit-1-\eps)}}\, dy. \eeq
We proceed as in Step \ref{sec:app}.3 and get that $|u(x)|\leq C
d(x,\partial\Omega)^{\theta}$ for some $\theta>\theta_0$. The claim
follows by induction.\hfill$\Box$

\medskip\noindent{\bf Step \ref{sec:app}.4:} We claim that $u\in
C^1(\overline{\Omega})$.

\smallskip\noindent{\it Proof of the claim:} With inequality
(\ref{ineq:claim3}), and the method used in Step \ref{sec:app}.2, we get
that
$$\lim_{k\to +\infty}u_k(x)=\int_\Omega
G(x,y)\left(\frac{|u|^{\crit-2-\eps}u(y)}{|\pi(y)|^s}-au\right)\, dy
$$
and
$$\lim_{k\to +\infty}\nabla u_k(x)=\int_\Omega
\nabla_xG(x,y)\left(\frac{|u|^{\crit-2-\eps}u(y)}{|\pi(y)|^s}-au\right)\, dy
$$
uniformly for $x\in\overline{\Omega}$. Since $u_k\rightharpoonup u$ in
$\huno$ when $k\to +\infty$, we get that $u\in
C^1(\overline{\Omega})$.\hfill$\Box$

\end{document}